\documentclass[12pt]{article}
\usepackage{epsfig,amssymb}

\makeatletter
\def\eqalign#1{\null\vcenter{\def\\{\cr}\openup\jot\m@th
  \ialign{\strut$\displaystyle{##}$\hfil&$\displaystyle{{}##}$\hfil
      \crcr#1\crcr}}\,}
\makeatother

\newcommand{\be}{\begin{equation}} 
\newcommand{\ee}{\end{equation}}
\newcommand{\beqa}{\begin{eqnarray}}
\newcommand{\eeqa}{\end{eqnarray}}
\newcommand{\bt}{\begin{theorem}}
\newcommand{\et}{\end{theorem}}
\newcommand{\bl}{\begin{lemma}}
\newcommand{\el}{\end{lemma}}
\newcommand{\bc}{\begin{corollary}}
\newcommand{\ec}{\end{corollary}}
\newcommand{\ba}{\begin{array}}
\newcommand{\ea}{\end{array}}

\newcommand{\la}{\label}
\newcommand{\ci}{\cite}

\newcommand{\cD}{{\cal D}}

\newcommand{\bi}{\bibitem}

\def\bbr{\mathbb R}
\def\complex{\mathbb C}
\def\bbc{\mathbb C}

\newtheorem{theorem}{Theorem}

\newtheorem{corollary}{Corollary}
\newtheorem{lemma}{Lemma}

\newcommand{\ti}{\tilde}
\newcommand{\wt}{\widetilde}
\newcommand{\de}{\delta}
\newcommand{\De}{\Delta}
\newcommand{\al}{\alpha}
\newcommand{\ga}{\gamma}

\newcommand{\si}{\sigma}
\newcommand{\Si}{\Sigma}
\newcommand{\om}{\omega}

\newcommand{\lb}{\lambda}
\newcommand{\ze}{\zeta}
\newcommand{\ka}{\varkappa}
\renewcommand{\th}{\theta}

\newcommand{\ep}{\varepsilon }

\newcommand{\Ai}{{\mathrm{Ai}}\,}

\begin{document}
\bigskip\bigskip\bigskip
\begin{center}
{\Large\bf 
Asymptotics of the Airy-kernel determinant
}\\
\bigskip\bigskip\bigskip

\centerline{ P.~ Deift}
\centerline{{\it Courant Institute of Mathematical Sciences,}}
\centerline{{\it New York, NY 10003, USA}}
\vskip .2in
\centerline{ A.~Its}
\centerline{{\it Department of Mathematical Sciences,}}
\centerline{{\it Indiana University -- Purdue University  Indianapolis}}
\centerline{{\it Indianapolis, IN 46202-3216, USA}}
\vskip .2in
\centerline{ I.~Krasovsky}
\centerline{{\it Department of Mathematical Sciences,}}
\centerline{{\it Indiana University -- Purdue University  Indianapolis}}
\centerline{{\it Indianapolis, IN 46202-3216, USA}}
\centerline{{\it and}}
\centerline{{\it Department of Mathematical Sciences}}
\centerline{{\it Brunel University}}
\centerline{{\it Uxbridge UB83PH}}
\centerline{{\it United Kingdom}}

\bigskip 
\bigskip\bigskip
\end{center}
\bigskip\bigskip\bigskip
\noindent{\bf Abstract.}
The authors use Riemann-Hilbert methods to compute the constant that arises 
in the asymptotic behavior of the Airy-kernel determinant of random matrix theory.

\newpage

\section{Introduction}
Let $K_s$ be the trace-class operator with kernel
\be
K_s(t,u) ={\Ai(t)\Ai'(u)-\Ai(u)\Ai'(t)\over t-u}
\ee
(see \ci{TW})
acting on $L^2(-s,\infty)$. Here $\Ai(x)$ is the Airy function (see, e.g., \ci{Abr}).
In this paper we are concerned with the behavior of $\det(I-K_s)$ as $s\to+\infty$.
Our main result is the following.

\noindent
{\bf Theorem 1.}
{\it
The large-$s$ asymptotic behavior of the Fredholm determinant $\det(I-K_s)$
is given by the formula
\be
\ln\det(I-K_s)=-{s^3\over 12}-{1\over 8}\ln s +\chi +O(s^{-3/2}),\la{Theorem1}
\ee
where
\be
\chi={1\over 24}\ln 2 +\ze'(-1),\la{const}
\ee
and $\ze(s)$ is the Riemann zeta-function. 
}

The Airy-kernel determinant $\det(I-K_s)$ is the edge scaling limit for the 
largest eigenvalue of a random $n\times n$ Hermitian matrix $H$ from 
the Gaussian Unitary 
Ensemble (GUE) (see \ci{Mehta,TW}) as $n\to\infty$: More precisely, 
if $\lb_1(H)\ge\lb_2(H)\ge\cdots\ge\lb_n(H)$ denote the eigenvalues of $H$,
then
\be\la{prob1}
\det(I-K_s)=\lim_{n\to\infty}\mbox{Prob }\{H\in \mbox{GUE}: 
(\lb_1(H)-\sqrt{2n})2^{1/2}n^{1/6}\le -s\}
\ee
(See \ci{F}, \ci{TW}, and also \ci{DG} for some history of (\ref{prob1})).

This determinant also describes the distribution of the longest increasing 
subsequence of random permutations \ci{BDJ,J}. Namely, let $\pi=i_1 i_2\cdots i_n$
be a permutation in the group $S_n$ of permutations of $1,2,\dots,n$. Then a subsequence
$i_{k_1},i_{k_2},\dots i_{k_r}$, $k_1<k_2<\cdots<k_r$,  
of $\pi$ is called an increasing subsequence of length
$r$ if $i_{k_1}<i_{k_2}<\cdots<i_{k_r}$. Let $l_n(\pi)$ denote the length of a longest
increasing subsequence of $\pi$ and let $S_n$ have the uniform probability distribution.
Then $l_n(\pi)$ is a random variable, and \ci{BDJ}
\be\la{prob2}
\det(I-K_s)=\lim_{n\to\infty}\mbox{Prob }\{\pi\in S_n: 
(l_n(\pi)-2\sqrt{n})n^{-1/6}\le -s\}
\ee

The distribution $F_{TW}(x)\equiv\det(I-K_{-x})$, 
known as the Tracy-Widom distribution,
admits the following integral representation \ci{TW}:
\begin{equation}\label{TW}
F_{TW}(x) = \exp \left \{ -\int_{x}^{\infty}(y-x)u^{2}(y)dy\right \},
\end{equation}
where $u(y)$ is the (global) Hastings-McLeod solution of the Painlev\'e II
equation
\begin{equation}\label{ds:30}
u''(y)=yu(y)+2u^3(y)\,,
\end{equation}
specified by the following 
asymptotic condition:
\begin{equation}\label{ds:40}
u(y)\sim \Ai(y)\qquad \mbox{as}\quad y\to+\infty.
\end{equation}
The behavior of $u(y)$ as $y \to -\infty$ is given
by the relation \cite{HM}:
\begin{equation}\label{ds:400}
u(y) = \sqrt{-\frac{y}{2}}
\left(1 + \frac{1}{8y^{3}} + O\left(y^{-6}\right)\right),
\qquad y \to -\infty,
\end{equation}
from which one learns that as $s\to +\infty$,
\be\la{asb}
-\int_{-s}^{\infty}(s+y)u^{2}(y)dy + 
{s^3\over 12}+{1\over 8}\ln s=as+b+o(1)
\ee
for some constants $a,b$. The content of Theorem 1 is that
$a=0$ and $b=\chi$ as in (\ref{const}). The value (\ref{const}) of the constant 
$\chi$ was conjectured by Tracy and Widom in \ci{TW} on the basis of the numerical 
evaluation of the l.h.s. of (\ref{asb}) as $s\to+\infty$ and by taking into account
the Dyson formula for a similar constant in the asymptotics of the so-called 
sine-kernel determinant \ci{TW}. The sine-kernel determinant describes the gap 
probability for GUE in the bulk scaling limit as $n\to\infty$ \ci{Mehta}.

Dyson's conjecture for the constant in the asymptotics of the 
sine-kernel determinant was proved rigorously in independent work by Ehrhardt \ci{Ehr}
and one of the authors \ci{K}, and a third proof was given later in \ci{DIKZ}. 
The two latter works use a Riemann-Hilbert-problem approach.
The proof in \ci{K} relies on a priori information from \ci{Widomarc},
whereas the proof in \ci{DIKZ} is self-contained. The proof of Theorem 1 in this paper
follows the method in \ci{DIKZ}.

As discussed in \ci{DIKZ}, the key difficulty in evaluating constants such as $\chi$
in (\ref{Theorem1}) in the asymptotic expansion of the determinants, is that 
in the course of the analysis one most 
naturally obtains expressions only for the logarithmic derivative with respect 
to some auxiliary parameter, say $\al$, in the problem, and not the determinant
itself. After evaluation of these 
expressions asymptotically, the constant of integration remains undetermined.
In \ci{DIKZ} and \ci{K}, this difficulty is overcome by utilizing a scaling limit of
finite-$n$ random matrices together with
universality in the sense of random matrix theory (see, e.g., \ci{Duniform}), in a way
that is inspired by, but different from, Dyson \ci{Dyson}.
We proceed as follows.

Consider the scaled Laguerre polynomials $p_k(x)$ 
defined for some integer $n$ by the orthogonality relation
\be
\int_0^\infty e^{-4nx}p_k(x)p_m(x)dx=\de_{k,m},\qquad k,m=0,1\dots,
\ee
The polynomial $p_k(x)=\ka_k x^k+\cdots$ is of degree $k$ and is related to the 
standard Laguerre polynomial $L^{(0)}_k(x)$ (see \ci{Szego}) as follows:
\[
p_k(x)=2\sqrt{n} L^{(0)}_k(4nx)
\]
with leading coefficient 
\be\label{ka-intr}
\ka_k=(-1)^k{2\sqrt{n}\over k!}(4n)^k\la{ka}.
\ee
The scaling here is chosen so that the asymptotic density of zeros of the 
polynomial $p_n(x)$ (with index $n$) is supported on the interval $(0,1)$
(as opposed to $(0,4n)$ for $L^{(0)}_n(x)$). See \ci{Szego,Dstrong} and below.

In the unitary random matrix ensemble defined by the Laguerre weight, the 
distribution function of the eigenvalues is given by the expression:
\be
dP(x_0,\dots,x_{n-1})={1\over C_n n!}\prod_{0\le i<j\le n-1}(x_i-x_j)^2
\prod_{j=0}^{n-1}e^{-4x_jn}dx_j,
\ee
where the normalization constant
\be\la{C}
C_n={1\over n!}\int_0^\infty\cdots\int_0^\infty
\prod_{0\le i<j\le n-1}(x_i-x_j)^2
\prod_{j=0}^{n-1}e^{-4x_jn}dx_j.
\ee
By a well known identity (e.g. \ci{Szego, Dbook}), the r.h.s. of the 
above expression gives
\be\la{C2}
C_n=\prod_{k=0}^{n-1}\ka_k^{-2}=(4n)^{-n^2}\prod_{k=0}^{n-1}k!^2,
\ee
where (\ref{ka}) was used.

For $\al\ge 0$,
the probability $D_n(\al)$ that the interval $(\al,\infty)$ has no
eigenvalues is given by
\be
D_n(\al)={1\over C_n n!}\int_0^\al\cdots\int_0^\al
\prod_{0\le i<j\le n-1}(x_i-x_j)^2
\prod_{j=0}^{n-1}e^{-4x_jn}dx_j,\la{Dint}
\ee
By standard arguments (cf. \ci{Dbook, Mehta}), this quantity can be written as
the Fredholm determinant of an integral operator on $L_2(0,\infty)$ in the following way:
\be
D_n(\al)=\det(I-K_n\chi_{(\al,\infty)}),\qquad
K_n(x,y)=\frac{1}{4}{\om_n(x)\om_{n-1}(y)-\om_n(y)\om_{n-1}(x)\over y-x},\la{det}
\ee
where
\be
\om_k(x)=e^{-2nx}p_k(x),\qquad k=0,1,\dots,
\ee
and $\chi_{(\al,\infty)}$ is the characteristic function of the interval 
$(\al,\infty)$.

If $x=1+1/(2n)+u/(2n)^{2/3}$ with $u$ fixed, then as $n\to\infty$, we obtain
from classical results on the asymptotics of the Laguerre polynomials 
(see \ci{Vanless,Szego}):
\be\eqalign{
\om_n\left(1+{1\over 2n}+{u\over (2n)^{2/3}}\right)=
\om_n\left({1\over 4n}(4n+2+2(2n)^{1/3}u)\right)=
(-1)^n{2\sqrt{n}\over (2n)^{1/3}}\left\{\Ai(u)+O(n^{-2/3})\right\};\\
\om_{n-1}\left(1+{1\over 2n}+{u\over (2n)^{2/3}}\right)=
\om_{n-1}\left({1\over 4n}\left[4(n-1)+2+2(2[n-1])^{1/3}\left(u+{2\over (2n)^{1/3}}
\right)+\right.\right.\\
\left.\left.O(n^{-2/3})\right]\right)=
(-1)^{n-1}{2\sqrt{n}\over (2[n-1])^{1/3}}\left\{\Ai\left(u+{2\over (2n)^{1/3}}\right)+
O(n^{-2/3})\right\},}\label{march17a}
\ee
where $\Ai(x)$ is the standard Airy function.
Let 
$$
K_{airy}(u,v) = 
{\Ai(u)\Ai'(v)-\Ai(v)\Ai'(u)\over u-v}.
$$
Set
$$
u^{(n)} = 1  + \frac{u}{(2n)^{2/3}} + \frac{1}{2n},\quad
v^{(n)} = 1  + \frac{v}{(2n)^{2/3}} + \frac{1}{2n}.
$$
It follows from (\ref{march17a}) that for any fixed $u$, $v$ we have
$$
\lim_{n \to \infty}\frac{1}{(2n)^{2/3}}K_{n}\left(u^{(n)}, v^{(n)}\right) 
=K_{airy}(u,v).
$$
In fact, this asymptotics is uniform for $u, v \geq L_{0}$,
where $L_{0}$ is an arbitrary constant. Indeed, for any $L_{0}$ there
exists $C = C(L_{0}) > 0$, $c = c(L_{0}) > 0$ such that
\begin{equation}\label{march17b}
\left|\partial^{j}_{u}\partial^{k}_{v}\left(
\frac{1}{(2n)^{2/3}}K_{n}\left(u^{(n)}, v^{(n)}\right)  - K_{airy}(u,v)\right)
\right|
\leq C\frac{e^{-cu}e^{-cv}}{n^{2/3}},
\end{equation}
$$
u, v \geq L_{0}, \quad j,k = 0, 1.
$$
This estimate can be proved in a same manner as estimate (3.8) in \cite{DG}.
In \cite{DG} the authors use global estimates for orthogonal polynomials
on ${\mathbb R}$ taken from \cite{Dstrong}: Here the relevant global estimates can be
obtained from \cite{Vanless}.

As in \cite{DG}, estimate (\ref{march17b}) immediately implies that
for any fixed $s\in{\mathbb R}$,
\be
\lim_{n \to \infty}D_{n}\left(1 - \frac{s}{(2n)^{2/3}}\right)
=\det\left(I-K_{airy}\chi_{(-s,\infty)})\right).\la{detA}
\ee

Below we obtain the asymptotics of the  determinant
$\det\left(I-K_{airy}\chi_{(-s,\infty)}\right) \equiv 
\det\left(I-K_{s}\right)$ as $s\to +\infty$.
In order to do this, we analyze the asymptotics of (\ref{det}) for all $\al$
from $\al$ close to zero to $\al=1-s/(2n)^{2/3}$.
Note that the determinant (\ref{det}) has the structure of so-called 
integrable determinants
\ci{IIKS}. Therefore, it is not surprising that there exists a differential identity
for ${d\over d\al}\ln D_n(\al)$ in terms of the solution of 
a related Riemann-Hilbert problem.
Solving the Riemann-Hilbert problem asymptotically as $n\to\infty$, 
we find the asymptotics of
this logarithmic derivative uniform for $\al\in[0,1-s_0/(2n)^{2/3}]$, $(2n)^{2/3}>s_0$
for some (large) $s_0>0$. 
Integrating these asymptotics from $\al$ close to zero to 
$\al=1-s/(2n)^{2/3}$, $s_0<s<(2n)^{2/3}$
we obtain the asymptotics of $D_n(1-s/(2n)^{2/3})$ 
provided we know the asymptotics of $D_n(\al)$ for $\al$ close to zero.
The latter, however, is readily obtained from the series expansion of the multiple 
integral formula for $D_n(\al)$ (see (\ref{Dint2},\ref{Dal0}) below).
More precisely, the ``inner workings'' of the method in this paper (cf. also
(133) in \ci{DIKZ}) can be seen from formula (\ref{intD}) below, which is obtained by
integrating the derivative $(d/d\al')\ln D_n(\al')$ from $\al'=\al_0$ to $\al'=\al$.
The key fact is that the estimate on the derivative is uniform for
$0\le\al'\le 1-s/(2n)^{2/3}$, $s>s_0$ (see (\ref{di1},\ref{di2})):
This leads to the error estimate $O(1/(n(1-\al)^{3/2}))$ in (\ref{intD}).
Using (\ref{Dal0}), we can then let $\al_0\to 0$: The singularities on the l.h.s. and 
the r.h.s. of (\ref{intD}) cancel out, and 
we are left with (\ref{intD2}). Using (\ref{intD2}), 
we immediately obtain Theorem 1. Note that in our calculations
formula (\ref{const}) for $\chi$ does not
arise from an evaluation of $D_n(\al_0)$ as $n\to\infty$ for some fixed $\al_0$. Rather
it arises, somewhat paradoxically, from the behavior of $D_n(\al_0)$ as $\al_0\to 0$
with $n$ {\it fixed} as given in (\ref{Dal0}).

In Section 2
the series expansion of $D_n(\al)$ for $n$ fixed and $\al\to 0$ is 
derived, as indicated above.
In Section 3 we obtain an asymptotic ($n\to\infty$) solution
of the Riemann-Hilbert problem related to (\ref{det}). Moreover, in Section 3, a differential 
identity for ${d\over d\al}\ln D_n(\al)$ is obtained in terms of the matrix elements 
(and their first derivatives) of the solution to the Riemann-Hilbert problem at the 
point $\al$. An alternative derivation of this identity, which is closer to the 
spirit of integrable systems and $\tau$-functions (see, e.g., \cite{BD,DIZ,IIKS}),
is given in the Appendix.
The identity is then evaluated asymptotically in Section 4 using asymptotics 
found in Section 3. In Section 5 the identity is integrated, and the results of Section 2
are then used to complete the proof.

\noindent {\bf Remark.}  Universality allows for considerable freedom in the choice
of the approximating ensemble in the above method. We choose to consider the Laguerre 
ensemble, although we could have considered, for example,
GUE itself: for GUE, however, the analysis turns out to be algebraically
more complicated. (For example, in the GUE case there will be two endpoints instead 
of one endpoint at $x=1$, see (\ref{ka-intr}) et seq.)
 In choosing the approximating ensemble, it is essential that 
the various constants that arise can be evaluated explicitly as in (\ref{Dal0}) 
and also in formula (17) in \ci{DIKZ}. In both cases we see that ultimately the 
formula for the
desired constant arises from classical formulae for the Legendre polynomials.   

In physics, and also in mathematical physics, universality is often viewed as 
a passive statement that certain systems ``behave in a similar fashion''. 
The thrust of this paper,
going back to Dyson \ci{Dyson}, is that universality can be used as an active 
analytical tool to obtain estimates for asymptotic problems of mathematical and 
physical interest.

\noindent {\bf Addendum.} We draw the attention of the reader to the work \cite{BBD}
of Baik, Buckingham, and DiFranco, in which the authors give a different proof of 
(\ref{const}) together with related results for GOE and GSE. The paper \cite{BBD}
appeared after our paper was written and refereed.

\section{Expansion of $D_n(\al)$ as $\al\to 0$.}
In this section we derive a series expansion for $D_n(\al)$ as $\al\to 0$. 
Changing the variables $x_j=(\al/2)(t_j+1)$ and expanding the exponent in (\ref{Dint}),
we obtain for fixed $n$:
\be\la{Dint2}
\eqalign{
D_n(\al)={1\over C_n n!}\left({\al\over 2}\right)^{n(n-1)+n}
\int_{-1}^1\cdots\int_{-1}^1
\prod_{0\le i<j\le n-1}(t_i-t_j)^2
\prod_{j=0}^{n-1}(1-2\al n(t_j+1)+O(\al^2))dt_j=\\
{1\over C_n}\left({\al\over 2}\right)^{n^2}A_n(1+O_n(\al)),}
\ee
where
\be
A_n={1\over n!}\int_{-1}^1\cdots\int_{-1}^1
\prod_{0\le i<j\le n-1}(t_i-t_j)^2\prod_{j=0}^{n-1}dt_j
\ee
can be expressed in terms of the product of the leading coefficients 
(cf. (\ref{C},\ref{C2})) of the Legendre polynomials:
\be
A_n=\prod_{k=0}^{n-1}
\frac{2^{2k}(k!)^4}{[(2k)!]^2}\frac{2}{2k+1}.
\ee
The asymptotics of $A_n$ as $n\to\infty$ (used first by Widom in \ci{Widomarc}, 
and then in \ci{DIKZ})
are given by the expression
\be\la{Aas}
\ln A_n=-n^2\ln 2 +n\ln(2\pi)  -{1\over 4}\ln n +   \frac{1}{12}\ln 2 + 3\zeta'(-1)
 + \wt\de_n,\quad n \to \infty.
\ee
where $\zeta'(x)$ is the derivative of Riemann's zeta-function, and 
$\wt\de_n\to 0$ as $n\to\infty$. The zeta-function originates 
from the expansion of the product of factorials.

The asymptotics of $C_n$ (\ref{C2}) have a similar form, 
\be\la{Cas}
\ln C_n= 
-\left({3\over 2}+\ln 4\right)n^2+ n\ln(2\pi)-{1\over 6}\ln n +2\ze'(-1)+
\hat\de_n,\qquad \hat\de_n\to 0,
\qquad n\to\infty.
\ee

Substituting the asymptotics (\ref{Aas},\ref{Cas}) into (\ref{Dint2}), we obtain
for $\al>0$:
\be\la{Dal0}
\ln D_n(\al)=\left({3\over 2}+\ln \al\right)n^2-{1\over 12}\ln{n\over 2}+
\ze'(-1)+\de_n+O_n(\al),
\ee
where $\de_n$ depends on $n$ only and $\de_n\to 0$ as $n\to\infty$. 
Note for later application (see proof of Lemma 2)
that the error term $O_n(\al)$ is analytic in $\al$, in particular,
$(d/d\al)O_n(\al)=O_n(1)$.  
We shall use formula (\ref{Dal0}) in the last section.

Caveat: $O_n(\al)\to 0$ as $\al\to 0$, $n$ fixed: no claim is made about $O_n(\al)$
as $n\to\infty$.

\section{Differential identity and the Riemann-Hilbert problem}
\subsection{Initial transformations}
In what follows, unless explicitly stated otherwise, we will always 
assume $0<\al< 1$.
At certain points in the text, however, we will also consider $\al$ in a small 
neighborhood $\cD_{\ep_0}(0)$ of $\al=0$ (see the discussion 
in the end of Section 3.1.)  

The multiple integral (\ref{Dint}) can be written as (cf. (\ref{C},\ref{C2})):
\be
D_n(\al)={1\over C_n}\prod_{j=0}^{n-1}\th_j^{-2},
\ee
where $\th_j$ are the leading coefficients of the polynomials
$q_j(x)=\th_j x^j+\cdots$ satisfying
\be\la{O}
\int_0^\al q_k(x)q_m(x)e^{-4nx}dx=\de_{km}, \qquad k,m=0,1,\dots
\ee

It is convenient to write this orthogonality relation in the form
\be\la{orth}
\int_0^\al q_j(x)x^k e^{-4nx}dx={\de_{jk}\over\th_j},\qquad 
k=0,1,\dots,j,\quad j=0,1,2,\dots
\ee
Note, in particular, that 
\be\la{aux}
\eqalign{
\int_0^\al q_j(x){\partial\over\partial\al}q_j(x)e^{-4nx}dx=\\
\int_0^\al
q_j(x)\left({d\th_j\over d\al}x^j+\mbox{polynomial of degree less than j}\right)
e^{-4nx}dx=
{1\over\th_j}{d\th_j\over d\al}.}
\ee

Using relation (\ref{aux}), we obtain
\be\la{D1}\eqalign{
{d\over d\al}\ln D_n(\al)={d\over d\al}\ln 
\prod_{j=0}^{n-1}\th_j^{-2}=
-2\sum_{j=0}^{n-1}{1\over\th_j}{d\th_j\over d\al}=
-2\sum_{j=0}^{n-1}
\int_0^\al q_j(x){\partial\over\partial\al}q_j(x)e^{-4nx}dx=\\
-\int_0^\al {\partial\over\partial\al}\left(\sum_{j=0}^{n-1}q^2_j(x)\right)
e^{-4nx}dx=
-{d\over d\al}\left(\int_0^\al \sum_{j=0}^{n-1}q^2_j(x)e^{-4nx}dx\right)+
\sum_{j=0}^{n-1}q^2_j(\al)e^{-4n\al}.}
\ee

By (\ref{O}) with $k=m=j$, the last integral (inside the brackets) in (\ref{D1})
equals $n$ and hence vanishes upon differentiation.

Applying the Christoffel-Darboux formula,
\be\la{CD}
\sum_{j=0}^{n-1}q^2_j(x)=
{\th_{n-1}\over \th_n}(q'_{n}(x)q_{n-1}(x)-q_n(x)q'_{n-1}(x)),
\ee
to the last sum in (\ref{D1}), we obtain
\be\la{idinterm}
{d\over d\al}\ln D_n(\al)={\th_{n-1}\over \th_n}e^{-4n\al}
(q'_{n}(\al)q_{n-1}(\al)-q_n(\al)q'_{n-1}(\al)).
\ee
Here and below the prime denotes differentiation w.r.t. the argument $x$.

Formula (\ref{idinterm}) shows that ${d\over d\al}\ln D_n(\al)$ depends only on $q_n$,
$q_{n-1}$. This property is crucial for the analysis below.

As noted in \ci{FIK}, orthogonal polynomials can be represented in terms of 
a solution to an associated Riemann-Hilbert problem. 
In the present case, the relevant Riemann-Hilbert problem is
formulated as follows:
Find a $2\times 2$ matrix-valued function $V(z)$ satisfying the conditions:
\begin{enumerate}
    \item[(a)]
        $V(z)$ is  analytic for $z\in\bbc \setminus [0,\al]$.
    \item[(b)] 
Let $x\in (0,\al)$. $V(z)$ has $L_2$ boundary values
$V_{+}(x)$ as $z$ approaches $x$ from above, and $V_{-}(x)$, from below.
They are related by the jump condition
\begin{equation}\la{Vjump}
            V_+(x) = V_-(x)
            \pmatrix{
                1 & e^{-4nx}\cr
                0 & 1},
            \qquad\mbox{$x\in (0,\al)$.}
        \end{equation}
    \item[(c)]
        $V(z)$ has the following asymptotic behavior as $z\to\infty$:
        \begin{equation}\la{Vinf} 
            V(z) = \left(I+ O \left( \frac{1}{z} \right)\right)z^{n\si_3},
            \qquad \mbox{where $\si_3=\pmatrix{1&0\cr 0&-1}$.}
    \end{equation}
\end{enumerate}
This Riemann-Hilbert problem (RHP) has a unique
solution for any $n$, $\al>0$, and, in particular, 
$V_{11}(z)=q_n(z)/\th_n$ and $V_{21}(z)=-2\pi i \th_{n-1} q_{n-1}(z)$.
Therefore we can rewrite the differential identity (\ref{idinterm}) in terms
of $V(z)$ in the form:\footnote{
An alternative derivation of
this identity is presented in the Appendix.}
\be\la{diffV}
{d\over d\al}\ln D_n(\al)={e^{-4n\al}\over 2\pi i}
(V_{11}(\al)V'_{21}(\al)-V'_{11}(\al)V_{21}(\al)).
\ee

In this section our task is to solve the RHP for $V(z)$
asymptotically (in other words, to find asymptotics of the polynomials
$q_k(z)$) as $n\to\infty$.
The results will then be used in Section 4.2. to evaluate the r.h.s. of (\ref{diffV}).

Following the steepest descent method for 
RH problems as described in \ci{Dstrong, Dbook},
we first of all need to find a 
so-called $g$-function: In the present situation this reduces to finding
a function analytic outside the interval $(-\infty,\al)$ 
and continuous up to the boundary with the properties: 
\begin{enumerate}
\item[(a)] $g(z)=\ln(z)+O(1/z)$ as $z\to\infty$;
\item[(b)] there exists a constant $l$ such that the boundary values $g_\pm(x)=\lim_{\ep\downarrow 0}
g(x\pm i\ep)$
of $g(z)$ are related as follows:
\be
g_+(x)+g_-(x)-4x-l=0,\qquad x\in (0,\al);\la{gpm}
\ee
\item[(c)] on $(0,\al)$, $g_+(x)-g_-(x)$ is purely imaginary, and $i(d/dx)(g_+(x)-g_-(x))>0$;
\item[(d)] $e^{g_+(x)-g_-(x)}=1$ on $(-\infty,0)$.
\end{enumerate}
A standard computation shows that if such a function $g(z)$ exists than it is unique.\footnote{
Note that as the contour for the RHP is $(0,\al)$,
the extra condition (4.14) for $g(z)$ in \ci{Dstrong} is redundant in the present situation.}

Formally, the derivative $g'(z)$ of $g(z)$  must have the properties:
\begin{enumerate}
\item[(a')] $g'(z)=1/z+O(1/z^2)$ as $z\to\infty$;
\item[(b')] $g_+'(x)+g_-'(x)=4$ for $x\in (0,\al)$.
\end{enumerate}

It is easy to verify that the following function satisfies these conditions:
\be\la{gprime}
g'(z)=2+{1+\al-2z \over \sqrt{z(z-\al)}}.
\ee
(In fact, $g'(z)$ is the unique function with $L^p$ boundary values $g_\pm '$ satisfying
(a') and (b') for any $1<p<2$.) In (\ref{gprime}), the branch is chosen so that 
$\sqrt{z(z-\al)}$ is analytic in the complement of $(0,\al)$ and positive for $z>\al$.

Therefore,
\[
g(z)=\int_\al^z g'(t)dt +C,
\]
where the constant $C$ is determined from the condition that $g(z)-\ln(z)=O(1/z)$
as $z\to\infty$. This gives
\be\la{g}
g(z)=2z-\al+\ln{\al\over 4}+\int_\al^z {1+\al-2t \over \sqrt{t(t-\al)}}dt,
\ee
and it is easy to verify that $g(z)$ indeed satisfies (a)--(d). 
From (\ref{gpm},\ref{g}) we now see that
\be\la{l}
l=-2\al+2\ln{\al\over 4}.
\ee

We need to analyze the RHP for $V(z)$ asymptotically as $n\to\infty$ uniformly for 
$0<\al<1-s_0/(2n)^{2/3}$ where $s_0$ is a fixed (large) number. The steepest descent method 
continues with the following steps (see \ci{Dstrong, Dbook}):

1) the RHP for $V$ is conjugated by $e^{ng(z)\si_3}$;

2) the contour $(0,\al)$ is split into lenses;

3) matching parametrices for the solution to the RHP are constructed (i) away from the end-points
$0$ and $\al$, (ii) in neighborhoods of $0$ and $\al$, respectively.

By means of these steps, the RHP reduces as 
$n\to\infty$ to a small norm problem which can be solved by a Neuman series.

All these steps go through in the standard way
except for the construction of the parametrix in a neighborhood of $\al$.
As we see from \ci{Dstrong,Dbook} the method requires that in a neighborhood $|z-\al|\le\ep$,
$\ep$ small and fixed,
\be\la{8+}
(g_+-g_-)(z)=(z-\al)^\beta (c+O(z-\al)),
\ee
for some $c\neq 0$ and some exponent $\beta>0$. (In \ci{Dstrong}, $\beta=3/2$.)
In our case for $0<z<\al$,
\be\la{8++}
(g_+-g_-)(z)=2\int_\al^z {1+\al-2t \over \sqrt{t(t-\al)}}dt=
{4\over\sqrt{\al}}(z-\al)^{1/2}(1-\al+O(z-\al)).
\ee
For any fixed $0<\al<1$ we see that $(g_+-g_-)(z)$ satisfies (\ref{8+}). 
As $\al\to 1$, we have to make the neighborhood $|z-\al|<\ep$ smaller and smaller.
The constant $c$ in (\ref{8+}) depends then on $\al$, but that, in itself, is not an insurmountable problem.
The real problem is that, unlike the situation in \ci{K},
the parametrix away from the points $0,\al$ (see \ci{Dstrong, Dbook}) contains certain terms
of the form $(z/(z-\al))^{1/4}$ evaluated on $\{z:|z-\al|=\ep\}$,
and as a result is not uniformly bounded when $1-\al$, and hence $\ep$, approach zero.
At the same time, there is not enough decay in the other relevant quantities 
to compensate for this.
The problem can be circumvented, however, by introducing a transformation
of the $z$-plane that ``regularizes'' the RHP in a neighborhood of $z=\al$.
Namely, set
\begin{figure}
\centerline{\psfig{file=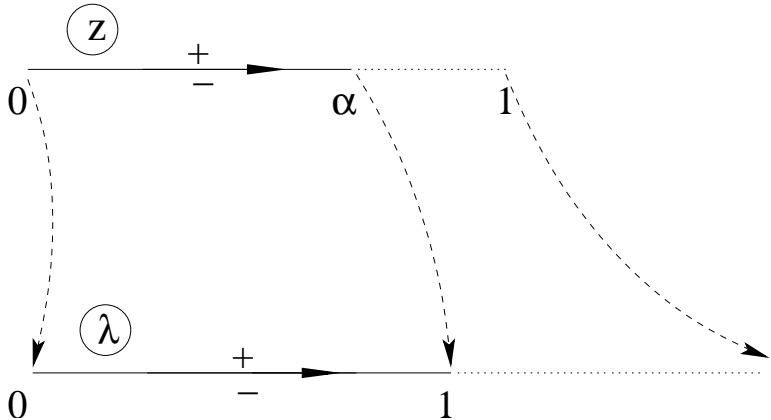,width=4.0in,angle=0}}
\vspace{0cm}
\caption{
Conformal mapping.}
\label{3}
\end{figure}
\be
\lb={1-\al\over\al}{z\over 1-z},\qquad z\neq1.
\ee
This fractional-linear transformation maps the interval $[0,\al]$ onto $[0,1]$,
the point $z=1$ is mapped to infinity, and infinity is mapped to $\lb=-(1-\al)/\al$.
The inverse transform is
\be\la{zlb}
z={\al\lb\over 1-\al+\al\lb},\qquad 
\lb\neq -{1-\al \over\al}.
\ee
Thus $z(\lb)$ is analytic from $\bbc\setminus\{ -(1-\al)/\al\}$
into $\bbc$, taking the complement of $[0,1]$
onto $\bbc\setminus [0,\al]$. 

The fact that in our case we could not obtain an estimate of the 
form (\ref{8+})
uniformly as $\al\uparrow 1$ originates in the vanishing of the numerator in the 
integral for $g_+-g_-$ in (\ref{8++}) at the point $t=(1+\al)/2\in (\al,1)$. 
Under the transformation 
$z\rightarrow\lb$ the point $(1+\al)/2$ is mapped to $\lb=1+\al^{-1}$. This point is at a
positive distance from the contour $0<\lb<1$ for $\al\in(0,1)$.
This means that we will be able to construct a parametrix for the solution of the RHP in 
the $\lb$ variable in a {\it fixed} neighborhood about $\lb=1$.
On the other hand, the point $\lb=-(1-\al)/\al$ (the image of $z$-infinity)
now approaches the contour as 
$\al\uparrow 1$, and we will need to contract the neighborhood of $\lb=0$ so that this 
point remains outside. We shall see, however, that this neighborhood presents no 
problem, as the relevant terms of the jump matrix for the final $R$-RHP 
(see (\ref{h00}) and the argument after (\ref{Rj}) below) decay 
sufficiently fast on the boundary of the neighborhood. 
 
For any $\lb\in \bbc\setminus([0,1]\cup\{-(1-\al)/\al\})$ set 
\be\la{U}
U(\lb)\equiv V(z(\lb)),
\ee 
where $z(\lb)=\al\lb/(1-\al+\al\lb)$ as in (\ref{zlb}). 

Then we obtain the following Riemann-Hilbert problem for $U(\lb)$:
%
%

\begin{enumerate}
    \item[(a)]
        $U(\lb)$ is  analytic for $\lb\in\overline{\bbc} \setminus 
                      ([0,1]\cup\{-(1-\al)/\al\})$.
    \item[(b)] 
Let $\lb\in(0,1)$.
$U$ has $L_2$ boundary values
$U_{+}(\lb)$ as $\lb$ approaches the real axis from
above, and $U_{-}(\lb)$, from below. 
They are related by the jump condition
\begin{equation}
            U_+(\lb) = U_-(\lb)
             \pmatrix{
                1 & e^{-4nz(\lb)}\cr
                0 & 1},
            \qquad\mbox{$\lb\in (0,1)$.}
        \end{equation}
    \item[(c)]
        $U(\lb)$ has the following asymptotic behavior as
$\lb\to -{1-\al\over\al}$ $(z\to\infty)$:
        \begin{equation} 
            U(\lb) = \left[I+ O \left( \frac{1}{z(\lb)} \right)\right]z(\lb)^{n\si_3}. 
    \end{equation}
\end{enumerate}           

We transfer $g(z)$ to the $\lb$-plane by defining
\be
\hat g(\lb)\equiv g(z(\lb)),\quad\mbox{for}\quad \lb\in\bbc\setminus
\left[-{1-\al \over\al},1\right].
\ee
Necessarily, $\hat g(\lb)$ is analytic on its 
domain. We obtain
\be
\hat g(\lb)=
2\al(\al-1){1-\lb\over 1-\al+\al\lb}+\al+\ln{\al\over 4}+
(1-\al)^{3/2}\int_{1}^\lb \frac{1+\al(1-t)}{(1-\al(1-t))^2}
\frac{dt}{\sqrt{t(t-1)}}.\la{49}
\ee
Note that $\hat g(\lb+0)-\hat g(\lb-0)=2\pi i$ on $(-(\al^{-1}-1),0)$ as this interval 
is the image of the half-axis $(-\infty,0)$ in the $z$-variable, 
where it is easy to conclude
(cf. (\ref{h(0)}) below) that $g_+(z)-g_-(z)=2\pi i$.
This jump in the $\lb$-variable is also
easy to obtain directly from (\ref{49}).

Let
\be\la{h}
h(\lb)=
{2(1-\al)^{3/2}\over e^{i\pi/2}}\int_{1}^\lb \frac{1+\al(1-t)}{(1-\al(1-t))^2}
\frac{dt}{\sqrt{t(1-t)}}
\ee
which is analytic in 
$\bbc\setminus ((-\infty,0)\cup(1,\infty))$. Here we choose the branch so that
$\sqrt{t(1-t)}$ is analytic in $\bbc\setminus((-\infty,0)\cup(1,\infty))$ and 
positive for $t\in (-1,1)$. The function $h(\lb)$ is the analytic continuation
of $\hat g(\lb+0)-\hat g(\lb-0)$ off the interval $(0,1)$.

Note that 
\be\la{h(0)}
h(0)={2\over e^{i\pi/2}}\int_1^0{1+\al-2\al x\over \sqrt{x(1-x)}}dx=2\pi i.
\ee

Now transform the RHP for $U$ as follows:
\be\la{T}
\wt T(\lb)=e^{-n l\si_3/2}U(\lb)e^{-n(\hat g(\lb)-l/2)\si_3},\qquad 
\lb\in\bbc\setminus\left([0,1]\cup\left\{-{1-\al\over\al}\right\}\right).
\ee
We easily obtain then that $\wt T(\lb)$ satisfies:
\begin{enumerate}
    \item[(a)]
        $\wt T(\lb)$ is  analytic for $\lb\in\overline{\bbc} \setminus [0,1]$.
    \item[(b)] 
For $\lb\in (0,1)$ the boundary values of $\wt T(\lb)$
are related by the jump condition
\begin{equation}
            \wt T_+(\lb) = \wt T_-(\lb)
            \pmatrix{
                e^{-n h(\lb)}& 1\cr
                0 & e^{n h(\lb)}},
            \qquad\mbox{$\lb\in (0,1)$.}
        \end{equation}
    \item[(c)]
        $\wt T(\lb)$ has the following asymptotic behavior as 
$\lb\to -{1-\al\over\al}$:
        \begin{equation} 
            \wt T(\lb) = I+ O \left(\lb+{1-\al\over\al}\right). 
    \end{equation}
\end{enumerate}
Note that the problem is now normalized to $I$ at $\lb=-{1-\al\over\al}$.

Since $\det\wt T(\lb)=1$ and $ \wt T(\lb)$ is analytic at infinity, it follows that
$ \wt T(\infty)$ is invertible.
The function $T(\lb)$ defined by 
\be
T(\lb)=\wt T(\infty)^{-1}\wt T(\lb)
\ee
is the solution to the same Riemann-Hilbert problem as $\wt T(\lb)$, 
with the $(c)$ condition replaced by 
\be
T(\lb)=I+O(1/\lb), \qquad \lb\to\infty.
\ee
Clearly,
\be\la{TT}
\wt T(\lb)=T^{-1}\left(-{1-\al\over\al}\right) T(\lb).
\ee

We now show that the RHP for $T$ is solvable for all $0\le\al<1$. 
For $0<\al<1$ the existence of such a $T(\lb)$ follows simply by
pushing forward $V(z)$, the solution of the RHP (\ref{Vjump},\ref{Vinf})
for the polynomials orthogonal on $(0,\al)$ with the weight $e^{-4nx}$:
the existence of $V(z)$ itself follows from the basic results of \ci{FIK,DZ3}.
So we are reduced to showing that $T(\lb)$ exists in the case $\al=0$ when 
the mapping $V(z)\to T(\lb)$ breaks down. 
For $\al=0$, $h(\lb)=4\ln(\sqrt{\lb}+i\sqrt{1-\lb})$, $0<\lb<1$. 
If $(\lb-1)^{1/2}$ (resp., $\lb^{1/2}$) denotes the branch which is analytic
in $\bbc\setminus [-\infty,1]$ (resp., $\bbc\setminus [-\infty,0]$),
then in particular $(\lb-1)^{1/2}_+=-(\lb-1)^{1/2}_-=i\sqrt{1-\lb}$, $0<\lb<1$,
and we find
\be
e^{n h(\lb)}=\left({(\lb-1)^{1/2}_+ +\lb^{1/2} \over (\lb-1)^{1/2}_- +\lb^{1/2}}
\right)^{2n}.
\ee
Thus if
$r(\lb)=((\lb-1)^{1/2}+\lb^{1/2})/2$, then
\be
\pmatrix{e^{-n h(\lb)}& 1\cr 0 & e^{n h(\lb)}}=
\pmatrix{(r_-/r_+)^{2n}& 1\cr 0& (r_+/r_-)^{2n}},\qquad 0<\lb<1.
\ee
Setting $Z(\lb)=T(\lb)r(\lb)^{2n\si_3}$, we see that $Z(\lb)$ solves the RHP: 
 \begin{enumerate}
    \item[(a)]
        $Z(\lb)$ is  analytic for $\lb\in\bbc\setminus [0,1]$.
    \item[(b)] 
For $\lb\in (0,1)$ the boundary values of $Z(\lb)$
are related by the jump condition
\begin{equation}
            Z_+(\lb) = Z_-(\lb)
            \pmatrix{
                1 & 1\cr
                0 & 1},
            \qquad\mbox{$\lb\in (0,1)$.}
        \end{equation}
    \item[(c)]
        $Z(\lb)$ has the following asymptotic behavior as 
$\lb\to\infty$:
        \begin{equation} 
            Z(\lb) = (I+ O(1/\lb))\lb^{n\si_3}. 
    \end{equation}
\end{enumerate}
This is the standard RHP for polynomials orthogonal on $(0,1)$ with the unit weight.
Therefore the desired solution $T(\lb)$ exists for $\al=0$ as well. 
This completes the proof of solvability of the RHP for $T(\lb)$ for all $0\le\al<1$.
The above proof of solvability for {\it all} $n$ is included only for completeness 
(cf. the last remark at the end of Section 3.5).

\begin{figure}
\centerline{\psfig{file=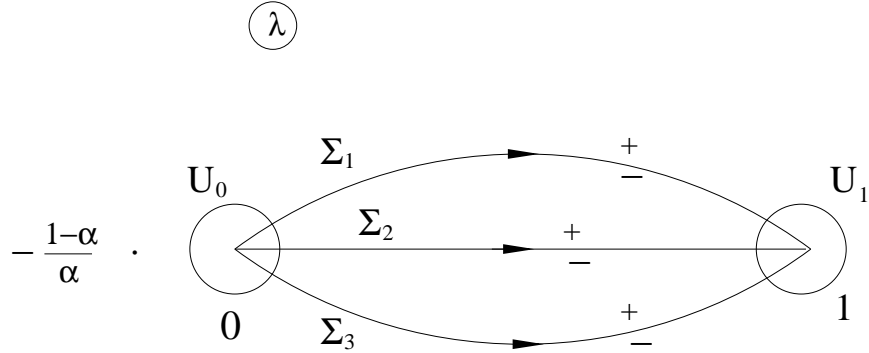,width=4.5in,angle=0}}
\vspace{0cm}
\caption{
Contour for the $S$-Riemann-Hilbert problem and the circular neighborhoods
$U_{1,0}$ of the points $1$, $0$. These neighborhoods will be introduced below 
in connection with the construction of parametrices.}
\label{5}
\end{figure}

As is standard in applications of the steepest descent method, we now deform the
RHP as follows. Let $\Si=\cup_{j=1}^3\Si_j$ be the oriented contour as in  Figure 
\ref{5}. Define a 
matrix-valued function $S(\lb)$ on $\bbc\setminus\Si$ by the expressions:
\be\la{S}
S(\lb)=
\cases{T(\lb),& for $\lb$ outside the lens,\cr
T(\lb)\pmatrix{1 & 0\cr -e^{-nh(\lb)}& 1}, & 
for $\lb$ in the upper part of the lens,\cr
 T(\lb)\pmatrix{1 & 0\cr e^{nh(\lb)}& 1},&
for $\lb$ in the lower part of the lens.
}
\ee

It is easy to verify that $S(\lb)$ solves the following RHP:

\begin{enumerate}
    \item[(a)]
        $S(\lb)$ is  analytic for $\lb\in\bbc \setminus\Sigma$, where 
$\Si=\cup_{j=1}^3\Si_j$.
    \item[(b)]  
The boundary values of $S(\lb)$ are related by the jump condition
\begin{equation}\la{Sjump}
\eqalign{   S_+(\lb) = S_-(\lb)
            \pmatrix{
                  1 & 0\cr
                  e^{\mp n h(\lb)} & 1},
            \qquad\mbox{$\lb\in(\Sigma_1\cup\Sigma_3)
\setminus\,$\{0,1\}},\\
\mbox{where the plus sign in the exponent is on 
$\Sigma_3$, and the minus sign, on $\Sigma_1$,}\\
            S_+(\lb) = S_-(\lb)
              \pmatrix{
                  0 & 1\cr
                  -1 & 0},
            \qquad\mbox{$\lb\in\Si_2\equiv(-1,1)$}.}
        \end{equation}
    \item[(c)]
$S(\lb)=I+O(1/\lb)$ as $\lb\to\infty$.
\end{enumerate}

For a fixed $0<\ep < 1/4$,  consider the circular neighborhood $U_1$  
of radius $\ep$ at the point 
$\lb=1$. Consider also the neighborhood $U_0$ 
of $\lb=0$ of radius $\ep_3(1-\al)$ 
for a fixed $1/2>\ep_3>0$. Note that $U_0$ contracts with growing $n$ for 
$\al=1-s_0/(2n)^{2/3}$. The point $-(1-\al)/\al$ lies outside $U_0$ for all
$\al\in(0,1)$. 

In $U_0$, we can expand the integrand in (\ref{h}) in powers of $t$ and $t/(1-\al)$:
\be\eqalign{
h(\lb)=h(0)+
{2\over e^{i\pi/2}\sqrt{1-\al}}\int_0^\lb
(1+\al-\al t)\left(1-{2\al t\over 1-\al}+O\left(t^2\over (1-\al)^2\right)\right)
\times\\
(1+t/2+O(t^2)){dt\over\sqrt{t}}=\\
h(0)+{4\sqrt{\lb}\over e^{i\pi/2}\sqrt{1-\al}}\left(1+\al+
{1-6\al-3\al^2\over 6(1-\al)}\lb+
O\left({\lb^2\over (1-\al)^2}\right)\right),\\
|\lb|\leq \ep_3(1-\al),}\la{h0}
\ee
uniformly in $\al$, and
where $h(0)=2\pi i$ (see (\ref{h(0)})).
It is the presence of $\sqrt{1-\al}$ in the denominator that will allow us to
construct a solution to the RHP using a contracting neighborhood $U_0$ as $\al$ 
approaches $1$. 

We shall now show that the the jump matrices for $S(\lb)$ 
on $\Si_1\cup\Si_3\setminus(U_1\cup U_0)$
are uniformly exponentially close to the identity (see (\ref{expest}) below)
as $n(1-\al)^{3/2}\to\infty$.

\begin{figure}
\centerline{\psfig{file=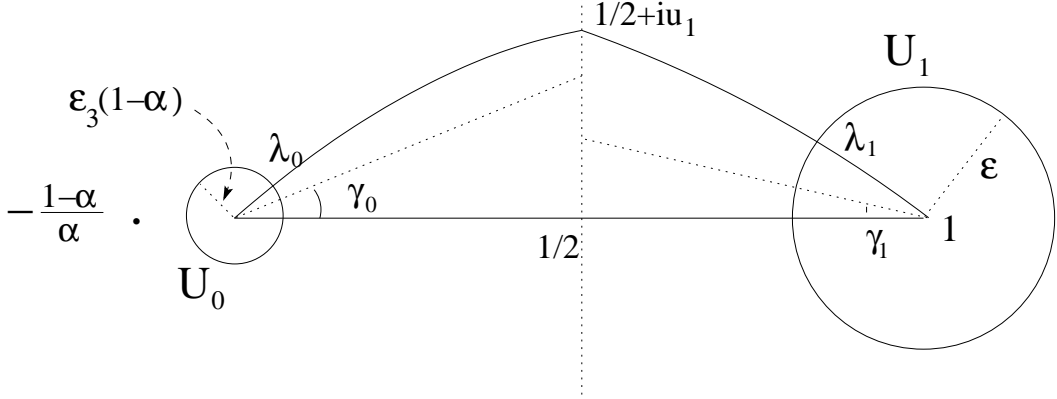,width=5.5in,angle=0}}
\vspace{0cm}
\caption{
Form of the contour for the $S$-Riemann-Hilbert problem.}
\label{6}
\end{figure}

To estimate the real part of $h(\lambda)$ outside of the neighborhoods 
$U_{0}$ and $U_1$, we now describe the form of the lens more precisely.
First, we assume that the contour $\Sigma_{3}$ is the
mirror image of $\Sigma_{1}$, i.e.
$$
\Sigma_{3} = \overline{\Sigma_{1}}.
$$ 
Therefore, we only need to describe the structure of
the contour $\Sigma_{1}$.
We assume that for $0\leq \Re \lambda \leq 1/2$ the contour $\Sigma_{1}$ 
lies above the straight line originating at zero,
and making a positive angle $\gamma_0$ with the real axis (see Figure \ref{6}).
The value of the angle $\gamma_0$  will be specified later on.
Similarly, the part of the contour between the vertical line
$\Re \lambda = 1/2$ and the boundary of the neighborhood 
$U_{1}$ lies above the line
$\Im \lambda = (1-\Re \lambda)\tan \gamma_{1}$ where, again,
the value of the angle $\gamma_{1}<\gamma_0$  will be specified later on.
Note that the contour $\Si$ has a well-defined limit as $\al\downarrow 0$.

Let $\lb_0$ (resp., $\lb_1$) be the point of intersection of the 
contour $\Sigma_{1}$ and the boundary of the disc $U_{0}$ (resp., $U_1$)
(see again Figure \ref{6}). 
Let $\Re \lb=\mu$, $\Im \lb=u$. Thus, $\lb=\mu+iu$, and on $\Si_1$,
$\Re\lb_0\le\mu\le\Re\lb_1$.
Fix some small $\ep_2>0$.
Suppose that $\Si_1$ and $\Si_3$ are so close to the real axis that
\be
{|u|\over\mu}<\ep_2,\qquad  {|u|\over 1-\mu}<\ep_2,\qquad \Re\lb_0\le\mu\le\Re\lb_1.
\ee
In particular, this implies that $\tan\ga_0<\ep_2$ and $\tan\ga_1<\ep_2$.
Furthermore, as $1+\al-\al\mu>1$, we have
\[
{|u|\over 1+\al-\al\mu}<|u|,
\]
and, as $1-\al+\al\mu>\al\mu$,
\[
{\al |u|\over 1-\al+\al\mu}<{|u|\over\mu}.
\]

The above inequalities  
allow us to perform the following estimate
on $h(\lb)$ for $\lb=\mu+iu$ in  $(\Si_1\cup\Si_3)\setminus \{U_1\cup U_0\}$. 
Using (\ref{h}), we obtain
\be\eqalign{
h(\lb)=h(\mu)+
{2(1-\al)^{3/2}\over e^{i\pi/2}}\int_{\mu}^{\mu+iu} 
\frac{1+\al-\al\mu-\al iv}{(1-\al+\al\mu+\al iv)^2}
\frac{d(\mu+iv)}{\sqrt{(\mu+iv)(1-\mu-iv)}}=\\
h(\mu)+{2(1-\al)^{3/2}\over\sqrt{\mu(1-\mu)}}
{1+\al-\al\mu\over (1-\al+\al\mu)^2}\int_0^u
\left(1-{i\al v\over 1+\al-\al\mu}\right)
\left(1+{i\al v\over 1-\al+\al\mu}\right)^{-2}\times\\
\left(1+{iv\over\mu}\right)^{-1/2}
\left(1-\frac{iv}{1-\mu}\right)^{-1/2}dv=
h(\mu)+{2(1-\al)^{3/2}\over\sqrt{\mu(1-\mu)}}
{1+\al-\al\mu\over (1-\al+\al\mu)^2}u
\left[1+O(\ep_2)\right],}\la{hout}
\ee
where the constant in the error term is uniform for $0\le\al<1$. 

The fraction $u/(1-\al+\al\mu)^2$ in the last 
equation of (\ref{hout}) can be estimated for some $\ep_4>0$ as 
\be\eqalign{
{|u|\over (1-\al+\al\mu)^2}> \ep\sin\ga_1>\ep_4,
\qquad \mbox{\rm for}\quad {1\over2}\le \mu\le \Re\lb_1,\\
{|u|\over (1-\al+\al\mu)^2}> {\mu\tan\ga_0\over(1-\al+\al\mu)^2}=
\frac{\tan\ga_0}{\mu(\al+(1-\al)/\mu)^2} >\\
\frac{\tan\ga_0}{(1+\ep_2/(\ep_3\sin\ga_0))^2}>\ep_4,
\qquad\mbox{\rm for}\quad  \Re\lb_0\le\mu\le{1\over 2},}\label{houtside}
\ee
where $\ep_4$ depends only on $\ep$ and $\ep_i$, $i=2,3$, $\ga_0$, $\ga_1$,
which in turn depend only on $\ep$, $\ep_2$, $\ep_3$.

Since $\Re h(\mu)=0$, we obtain
from (\ref{houtside}) as $n\to\infty$ for
sufficiently small $\ep_2>0$:
\be\eqalign{
|e^{-nh(\lb)}|=O(e^{-\rho c}),\qquad 
\lb\in\Si_1\setminus(U_0\cup U_1),\\
|e^{nh(\lb)}|=O(e^{-\rho c}),\qquad 
\lb\in\Si_3\setminus(U_0\cup U_1)}\la{expest}
\ee
{\it uniformly} for $\al\in[0,1-s_0/(2n)^{2/3}]$ for some (large) $s_0>0$ and 
all $n>s_0^{3/2}/2$, for some $c=c(\ep,\ep_2,\ep_3)>0$, where 
$$
\rho=n|1-\al|^{3/2}.
$$
So except for the jump on the interval $(0,1)$ and the jumps inside $U_{1}$, $U_0$,
the jumps of $S(\lb)$ are indeed exponentially close to the identity as $\rho\to\infty$.

For later purposes, we shall need the series expansion of $h(\lb)$ at $\lb=0,1$.
We have:
\be\la{h00}
h(\lb)=2\pi i+
{4\sqrt{\lb}\over e^{i\pi/2}\sqrt{1-\al}}
\left(1+\al+{1-6\al-3\al^2\over 6(1-\al)}\lb+
O\left({\lb^2\over (1-\al)^2}\right)\right),\qquad \lb\to 0;
\ee
\be\la{h1}\eqalign{
h(\lb)=4(1-\al)^{3/2}\sqrt{u}\left(1-(\al+1/6)u+
(\al^2+3\al/10+3/40)u^2+O(u^3)\right),\qquad
\lb=1+u,\\ u\to 0.}
\ee

In (\ref{h00}) the cut of the root lies to the left of $\lb=0$, and $-\pi<\arg\lb<\pi$,
whereas in (\ref{h1}) the cut lies to the right of $\lb=1$, and $0<\arg u< 2\pi$.

Note the crucial fact that, as follows from
(\ref{h00}),  (\ref{h1}), the quantity
$n|h(\lb)|$ (resp., $n|h(\lb)-2\pi i|$)  is {\it uniformly large}  on
the boundary 
$\partial U_1$ (resp., $\partial U_0$)  for 
some (large) $s_0>0$ for all $\al\in[0,1-s_0/(2n)^{2/3}]$,
if $(2n)^{2/3}>s_0$.
Indeed, it is of order $s_0^{3/2}$ for $\lb$ on $\partial U_1$
(resp., of order $n$ for $\lb$ on $\partial U_0$).
This will allow us to obtain the desired asymptotic solution of the Riemann-Hilbert 
problem.

For technical reasons (see the end of the section 4.2. below and also proof of 
Corollary 2 in \ci{DIKZ}), we need to control the solution of the RHP for all
$\alpha \in {\cal D}_{\ep_{0}}(0)
\cup [\ep_0,1-s_0/(2n)^{2/3}]$, where ${\cal D}_{\ep_{0}}$ denotes the disc
of radius $\ep_{0}$ about zero in the complex $\alpha$-plane with $\ep_0$ small.
For all  $\alpha \in {\cal D}_{\ep_{0}}(0)$ we use the fixed contour $\Si=\Si_{\al=0}$
in Figure 3 corresponding to $\al=0$. By the preceding calculation we see that
$|\Re h(\lb;\al=0)|\ge c_0>0$ for all 
$\lambda \in (\Si_1\cup\Si_3)\setminus (U_0\cup U_1)$.
Thus
\be 
|e^{-nh(\lb;\al=0)}|\le e^{-nc_0},\qquad 
\lambda \in (\Si_1\cup\Si_3)\setminus (U_0\cup U_1).
\ee
Hence, by continuity, we must have
\be
|e^{-nh(\lb,\al)}|\le e^{-nc'_0}
\ee
for all $\lambda \in (\Si_1\cup\Si_3)\setminus (U_0\cup U_1)$
and all $\alpha \in {\cal D}_{\ep_{0}}(0)$, $0<c'_0<c_0$, $\ep_0$ sufficiently small.

We now begin the construction of parametrices which give, in their respective regions,
the leading contribution to the asymptotics for the RHP.

\subsection{Parametrix in $\bbc\setminus(\overline{U_{1}\cup U_{0}})$}

First, because of the exponential convergence described above, we expect the 
following model problem to play a role in constructing a parametrix for the 
solution of the RHP as $n\to\infty$:
\begin{enumerate}
    \item[(a)]
        $N(\lb)$ is  analytic for $\lb\in\bbc \setminus[0,1]$,
    \item[(b)]
\be
N_{+}(\lb) = N_{-}(\lb)
            \pmatrix{0&1\cr
              -1&0},
\qquad\mbox{$\lb \in (0,1)$},
\ee
\item[(c)]
\be
N(\lb) = I+ O \left( \frac{1}{\lb} \right), 
     \qquad \mbox{as $\lb\to\infty$.}
\ee
\end{enumerate}
The solution $N(\lb)$ can be found in the standard way by first
transforming $N(\lb)$ with a $2\times 2$ unitary transformation to the form 
for which the jump matrix is diagonal and
then solving the two resulting scalar Riemann-Hilbert problems (cf. \ci{Dbook}).
We obtain
\be
N(\lb)={1\over 2}
\pmatrix{m+m^{-1}&-i(m-m^{-1})\cr i(m-m^{-1})& m+m^{-1}},
\qquad
m(\lb)=\left({\lb-1\over \lb}\right)^{1/4},\la{N}
\ee
where $m(\lb)$ is analytic outside $[0,1]$ and $m(\lb)\to +1$ as $\lb\to\infty$.
Note that $\det N(\lb)=1$ and that $N(\lb)$ is the unique $L^p$ solution of the RHP for
any $1<p<4$.


\subsection{Parametrix at $\lb=1$}

Now let us construct a parametrix in $U_1$.
We look for an analytic matrix-valued function $P_1(\lb)$ in $U_1$ which has
the same jump
relation as $S(\lb)$ on $\Si\cap U_1$ and instead of a condition at infinity
satisfies the matching condition on the boundary
\be\la{match2}
P_1(\lb)N^{-1}(\lb)=I+O(1/\rho), \qquad \lb\in\partial U_1, \qquad
\rho=n|1-\al|^{3/2},
\ee
uniformly in $\lb$ and $\al$ as $\rho\to\infty$.

Define:
\be
\phi(\lb)=\cases{e^{i\pi}h(\lb)/2,& for $\Im \lb>0$,\cr
h(\lb)/2,& for $\Im \lb<0$}.
\ee
This function is analytic in $U_1$ outside $(1-\ep,1]$.

We look for $P_1(z)$ in the form:
\be\la{P1}
P_1(\lb)=E_n(\lb)\hat P(\lb)e^{n\phi(\lb)\si_3},
\ee
where $E_n(\lb)$ is analytic and invertible ($\det E_n\neq 0$) in 
a neighborhood of $U_1$,
and therefore does not affect the jump and analyticity conditions for 
$\hat P(\lb)e^{n\phi(\lb)\si_3}$.

As $P_1(\lb)$ is required to satisfy the jump relations (\ref{Sjump}) for $S$,
it is easy to verify that $\hat P(\lb)=E_n(\lb)^{-1}P_1(\lb)e^{-n\phi(\lb)\si_3}$ 
satisfies jump conditions with {\it constant} jump matrices:
\be\eqalign{
\hat P_{+}(\lb) = \hat P_{-}(\lb)
              \pmatrix{
                  1 & 0\cr
                  1 & 1},
            \qquad\mbox{$\lb\in((\Si_1\cup\Si_3)\cap U_1)\setminus\{ 1\}$,}\\
\hat P_{+}(\lb) =\hat  P_{-}(\lb)
            \pmatrix{
                0 & 1\cr
                -1 & 0},\qquad\mbox{$\lb\in\Si_2\cap U_1$.}}
\ee

Now introduce a mapping of $U_1$ onto a new $\ze$-plane
\be\la{ze1}
\ze=n^2\phi(\lb)^2=4n^2(1-\al)^3 u (1-(2\al+1/3)u+
(3\al^2+14\al/15+8/45)u^2 + O(u^3)),\qquad \lb=1+u,
\ee
where we used (\ref{h1}). The expansion at $\lb=1$ is uniform for $\al$ in a bounded set.

Choosing a sufficiently small $\ep>0$, 
we see that $\zeta(\lb)$ is analytic and one-to-one in the
neighborhood $U_1$.

Note that if $\al\in[0,1-s_0/(2n)^{2/3}]$ then $|\ze|=O(\rho^2)$ uniformly 
large, if $s_0$ is large, 
on the boundary $\partial U_1$ and in $\al$. This is a crucial fact in the present work. 
When  $\al=1-s_0/(2n)^{2/3}$, we have $\rho=s_0^{3/2}/2$.

\begin{figure}
\centerline{\psfig{file=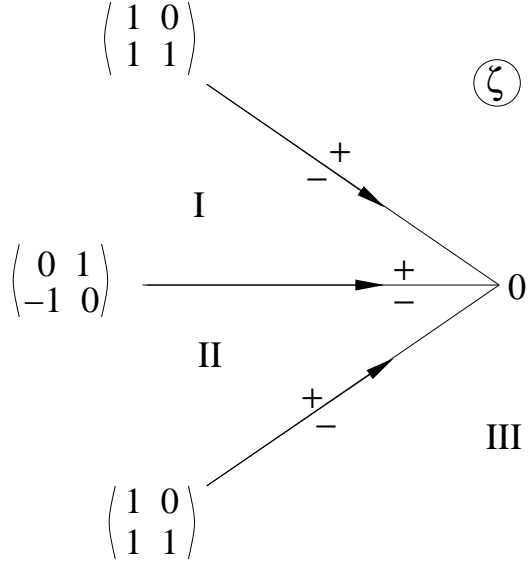,width=2.7in,angle=0}}
\vspace{0cm}
\caption{
Contour of the Riemann-Hilbert problem for $\Psi(\ze)$ (the case of $U_1$).}
\label{8}
\end{figure}

Let us now choose the exact form of the contours in $U_{1}$ so that their images
under the mapping $\zeta(\lb)$ are straight lines (see Figure \ref{8}).
Set
\be\la{PPsi1}
\hat P(\lb)=\Psi(\zeta),
\ee
So the jump matrices for $\Psi(\zeta)$ 
are the same as for $\hat P(\lb)$ (they are shown in Figure \ref{8}).
A matrix $\Psi(\ze)$  satisfying these jump conditions was constructed in \ci{KMVV}
in terms of Bessel functions, namely:

1) region I
\begin{equation}\label{Q1}
    \Psi(\zeta) ={1\over 2}
    \pmatrix{
        H_{0}^{(1)}(e^{-i\pi/2}\zeta^{1/2}) &
        H_{0}^{(2)}(e^{-i\pi/2}\zeta^{1/2}) \cr
        \pi \zeta^{1/2} \left(H_{0}^{(1)}\right)'(e^{-i\pi/2}\zeta^{1/2}) &
        \pi \zeta^{1/2} \left(H_{0}^{(2)}\right)'(e^{-i\pi/2}\zeta^{1/2})
    },
\end{equation}

2) region II
\begin{equation}\label{Q2}
    \Psi(\zeta) ={1\over 2}
    \pmatrix{
        H_{0}^{(2)}(e^{i\pi/2}\zeta^{1/2}) &
        -H_{0}^{(1)}(e^{i\pi/2}\zeta^{1/2}) \cr
        -\pi \zeta^{1/2} \left(H_{0}^{(2)}\right)'(e^{i\pi/2}\zeta^{1/2}) &
        \pi \zeta^{1/2} \left(H_{0}^{(1)}\right)'(e^{i\pi/2}\zeta^{1/2})
    },
\end{equation}

3) region III
\begin{equation}\label{Q3}
    \Psi(\zeta) =
    \pmatrix{
     I_{0} (\zeta^{1/2}) & \frac{i}{\pi} K_{0}(
        \zeta^{1/2})\cr
     \pi i \zeta^{1/2} I_{0}'(\zeta^{1/2}) & 
     -\zeta^{1/2} K_{0}'(\zeta^{1/2})},
\end{equation}
where $-\pi<\arg(\zeta)<\pi$.

Here the square root $\sqrt{\zeta}$ has the cut on $(-\infty,0)$. Hence,
$\sqrt{\zeta}=-n\phi(\lb)$ for $-\pi<\arg(\zeta)<\pi$.

The large-$\ze$ asymptotics of Bessel functions give 
(here we choose $s_0$, depending only on $\ep$,
sufficiently large): 
\be\eqalign{
\Psi(\ze)={1\over\sqrt{2}}(\pi\sqrt\ze)^{-\si_3/2}\pmatrix{1&i\cr i&1}
\left[I+{1\over 8\sqrt\ze}\pmatrix{-1&-2i\cr -2i&1}-\right.\\ \left.
{3\over 2^7\ze}\pmatrix{1&-4i\cr 4i&1}+
O(\ze^{-3/2})\right]e^{\sqrt\ze\si_3}}\la{asPsi1}
\ee
uniformly on the boundary $\partial U_1$.

Thus 
\be\la{P1p}
P_1(\lb)=E_n(\lb)\Psi(\ze(\lb))e^{n\phi(\lb)\si_3},
\ee
where the function $E_n(\lb)$ is found from the matching condition
to be
\be
E_n(\lb)=
{1\over\sqrt{2}}N(\lb)
\pmatrix{1&-i\cr-i&1}(\pi \sqrt{\ze})^{\si_3/2}.\la{E1}
\ee

Now to complete the construction of the parametrix it only remains to show that $E_n(\lb)$
is an analytic function in $U_1$ (clearly, $\det E_n(\lb)\neq 0$). 
First, we show that it has no jump on the real $\ze$-axis.
This is easy to verify using the jump condition for $N(\lb)$ and the identity 
$\ze_-=\ze_+ e^{-2\pi i}$ on the negative half axis.
Moreover, a simple calculation shows that $E_n(\lb)$ has no pole at $\lb=1$.
Thus, $E_n(\lb)$ is analytic in $U_1$, and the parametrix in $U_1$ is 
given by the equations (\ref{P1},\ref{PPsi1},\ref{Q1},\ref{Q2},\ref{Q3},\ref{E1})
for $\al\in [0,1)$.

Below we shall need the first three terms in the matching condition for $P_1$. 
Using (\ref{asPsi1}), we obtain
\be\la{j1}
P_1(\lb)N^{-1}(\lb)=I+\De_1(\lb)+\De_2(\lb)+O\left({1\over\rho^3}\right),
\qquad \lb\in\partial U_1.
\ee
Here
\be\eqalign{
\De_1(\lb)={1\over
  8\sqrt{\ze}}N(\lb)\pmatrix{-1&-2i\cr-2i&1}N(\lb)^{-1}=
{1\over 16\sqrt{\ze}}\pmatrix{-3m^2+m^{-2}&-i(3m^2+m^{-2})\cr
-i(3m^2+m^{-2})&3m^2-m^{-2}},\\
\De_2(\lb)=
{3\over 2^7\ze}N(\lb)\pmatrix{-1&4i\cr-4i&-1}N(\lb)^{-1}=
{3\over 2^7\ze}\pmatrix{-1&4i\cr-4i&-1},}\la{De1}
\ee
where $m(\lb)$ is defined in (\ref{N}).
Note that both $\De_1(\lb)$ and $\De_2(\lb)$ 
are meromorphic functions in $U_1$ with a simple pole at $\lb=1$.

Recall that we use the contour $\Si=\Si_{\al=0}$ for all $\al\in D_{\ep_0}(0)$,
$\ep_0$ small. For such $\al$, the map $\lb\to\ze$ maps $U_1$ (consisting
of the three regions separated by $\Si$) onto a set,
region to region, where the lines separating each region are now no longer 
straight but lie in small cones about the original ones. 
The opening angles of the cones are proportional to $|\Im\al|$. 
Using the same definition for $\Psi$ as in (\ref{Q1}--\ref{Q3}) 
for each of the new regions
I, II, III, we find again that (\ref{asPsi1}) is valid, and that
$P_1(\lb)N^{-1}(\lb)$ has the same expansion (\ref{j1}) as in the case $0\leq\al<1$.
Note that the values of $\ep_0$ and $s_0$ can be changed (now and below)
if necessary.


\subsection{Parametrix at $\lb=0$}

The construction of the parametrix in $U_{0}$ is similar. 
Recall, however, that the radius of $U_0$ is $\ep_3(1-\al)$, 
so it decreases as $\al\to 1$,
i.e. as the pole of $h(\lb)$ approaches the point $\lb=0$.
We shall see that this neighborhood produces asymptotics for the RHP
in inverse powers of $n$. 

We look for an analytic matrix-valued function $P_0(z)$ 
in the neighborhood $U_{0}$
which satisfies the same jump conditions as $S(\lb)$ on
$\Sigma\cap U_0$, and satisfies the matching condition
\be\la{match0}
P_0(\lb)N^{-1}(\lb)=I+(1-\al)^{-1/2}O(1/n)
\ee
uniformly in $\lb$ on the boundary $\partial U_0$ as $n\to\infty$.

Below we define functions in $U_0$ which play the same role 
as $\phi$, $E_n$, and $\hat P$
in $U_1$. We use the same notation for these quantities as before.  
Namely, let
\be
\phi(\lb)=\cases{e^{i\pi}(h(\lb)-2\pi i)/2,& for $\Im \lb>0$,\cr
(h(\lb)-2\pi i)/2,& for $\Im \lb<0$}.
\ee
This function is analytic in $U_{0}$ outside $[0,\ep_3(1-\al)]$.

As above, we look for the parametrix $P_0(\lb)$ in the form:
\be\la{P0}
P_0(\lb)=E_n(\lb)\hat P(\lb)e^{n\phi(\lb)\si_3},
\ee
We obtain that 
\be\eqalign{
\hat P_{+}(\lb) = \hat P_{-}(\lb)
              \pmatrix{
                  1 & 0\cr
                  1 & 1},
            \qquad\mbox{$\lb\in((\Si_1\cup\Si_3)\cap U_0)\setminus\{ 0\}$,}\\
\hat P_{+}(\lb) =\hat  P_{-}(\lb)
            \pmatrix{
                0 & 1\cr
                -1 & 0},\qquad\mbox{$\lb\in\Si_2\cap U_{0}$.}}
\ee

We choose the following mapping of $U_{0}$ onto a $\ze$-plane
\be\la{ze0}
\ze=n^2\phi(\lb)^2=e^{-i\pi}4n^2(1+\al)^2{\lb\over 1-\al}\left(1+
{1-6\al-3\al^2\over 3(1-\al^2)}\lb+
O\left({\lb^2\over (1-\al)^2}\right)\right),
\ee
where we used (\ref{h00}).

Choosing a sufficiently small $\ep_3>0$, 
we see that $\zeta(\lb)$ is analytic and one-to-one in the neighborhood $U_0$.

Let us also choose the exact form of the contours in $U_0$ so that their images
under the mapping $\zeta(\lb)$ are direct lines. 
In the $\ze$-plane the contour and the jump matrices for $\hat P(\lb)$ are the same as
in Figure \ref{8} with the only difference that all directions are reversed
(pointing away from $\ze=0$).
It is easily seen that the function
\be\la{PPsi0}
\hat P(\lb)=\si_3\Psi(\zeta)\si_3,
\ee
where $\Psi(\zeta)$ is given by (\ref{Q1}--\ref{Q3})
satisfies the jump conditions in this case.

Finally, we calculate $E_n$ and obtain
\be
P_0(\lb)=E_n(\lb)\si_3\Psi(\ze(\lb))\si_3 e^{n\phi(\lb)\si_3},
\ee
where 
\be
E_n(\lb)=
{1\over\sqrt{2}}N(\lb)
\pmatrix{1&i\cr i&1}(\pi \sqrt{\ze})^{\si_3/2}
\ee
(the analyticity of $E_n(\lb)$ in $U_0$ is verified as above).

Then we see immediately from (\ref{asPsi1},\ref{ze0}) that 
\be\la{match0a}
P_0(\lb)N^{-1}(\lb)=E_n(\lb)\si_3\Psi(\ze)\si_3 e^{n\phi(\lb)\si_3}N^{-1}(\lb)=
I+{1\over\sqrt{\lb}}O\left({1\over\sqrt{\ze}}\right)=
I+{1\over\sqrt{1-\al}}O\left({1\over n}\right)
\ee
uniformly in $\lb\in\partial U_0$ and $\al\in[0,1-s_0/(2n)^{2/3}]$.
Of course, the bound in (\ref{match0a}) blows up if
$\al\to 1$ too rapidly: for $0\le\al<1-s_0/(2n)^{2/3}$, we see that the error
term is $O(n^{-2/3})$.

Thus the construction of the parametrix in $U_0$ is now complete.

Using the expansion of $\Psi(\ze)$, we can extend (\ref{match0a}) to a full asymptotic 
series in inverse powers of $n$. Substituting (\ref{asPsi1}) into (\ref{match0a}), 
we obtain in particular:
\be\la{j0}
P_0(\lb)N^{-1}(\lb)=I+\De_1(\lb)+\De_2(\lb)+
{1\over\sqrt{1-\al}}O\left({1\over n^3}\right),
\ee
where
\be\eqalign{
\De_1(\lb)={1\over 8\sqrt{\ze}}N(\lb)\pmatrix{-1&2i\cr
  2i&1}N(\lb)^{-1}=
{1\over 16\sqrt{\ze}}\pmatrix{m^2-3m^{-2}&i(m^2+3m^{-2})\cr
i(m^2+3m^{-2})& -m^2+3m^{-2}},\\
\De_2(\lb)=
{3\over 2^7\ze}N(\lb)\pmatrix{-1&-4i\cr4i&-1}N(\lb)^{-1}=
{3\over 2^7\ze}\pmatrix{-1&-4i\cr4i&-1}.}\la{De0}
\ee
As above, note that $\De_1(\lb)$ and $\De_2(\lb)$ 
are meromorphic functions in $U_0$ with a simple pole at $\lb=0$.

For sufficiently small $\ep_0$, the estimate (\ref{j0}) extends uniformly for
$\al\in\cD_{\ep_0}(0)\cup [0,1-s_0/(2n)^{2/3}]$ for all $n>s_0^{3/2}/2$,
and $\lb\in\partial U_0$ as in Section 3.3.


\subsection{Final transformation of the problem}

\begin{figure}
\centerline{\psfig{file=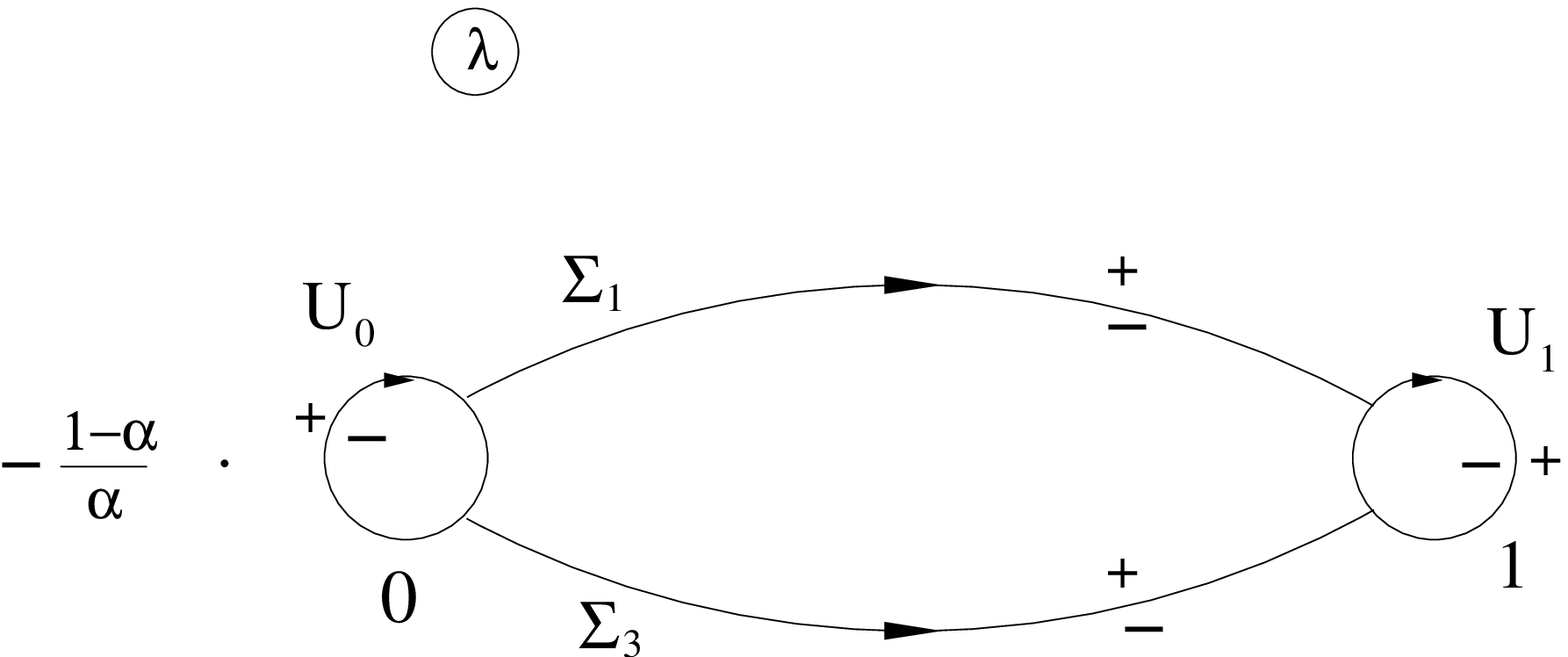,width=4.5in,angle=0}}
\vspace{0cm}
\caption{
Contour for the $R$-Riemann-Hilbert problem.}
\label{9}
\end{figure}

Now construction of the parametrices is complete, and we are ready
for the last transformation of the Riemann-Hilbert
problem. Let
\be
R(\lb)=\cases{S(\lb)N^{-1}(\lb),& 
$\lb\in \bbc\setminus (\overline{U_0\cup U_1}\cup\Si)$,\cr
S(\lb)P_0^{-1}(\lb),& 
$\lb\in U_0\setminus\Si$,\cr
S(\lb)P_1^{-1}(\lb),&
$\lb\in U_1\setminus\Si$.}\la{Rdef}
\ee
It is easy to see that this function has jumps only on 
$\partial U_1$, $\partial U_0$, and parts of $\Si_1$, $\Si_3$
lying outside the neighborhoods $\overline{U_1}$, $\overline{U_0}$ (we
denote these parts $\Si^\mathrm{out}$).
The contour is shown in Figure \ref{9}. Outside this contour, $R(\lb)$ is analytic.
Besides, $R(\lb)=I+O(1/\lb)$ as $\lb\to\infty$. 

The jumps are as follows:
\be\eqalign{
R_+(\lb)=R_-(\lb)N(\lb)\pmatrix{1&0\cr e^{\mp n h(\lb)}& 1}N(\lb)^{-1},
\qquad \lb\in\Si^\mathrm{out}_1\cup \Si^\mathrm{out}_3,\\
\mbox{where the ``-'' sign in the exponent is taken on $\Si^\mathrm{out}_1$, and ``+'',
on $\Si^\mathrm{out}_3$,}\\
R_+(\lb)=R_-(\lb)P_{0}(\lb)N(\lb)^{-1},\qquad \lb\in \partial U_0\setminus
\mbox{\{intersection points\}},\\
R_+(\lb)=R_-(\lb)P_1(\lb)N(\lb)^{-1},\qquad \lb\in \partial U_1\setminus
\mbox{\{intersection points\}}.}\la{Rj}
\ee

The jump matrix on $\Si^\mathrm{out}$ can be uniformly estimated
(both in $\lb$ and $\al\in[0,1-s_0/(2n)^{2/3}]$)
as $I+O(\exp(-c\rho))$, where $c$ is a positive constant.
In view of the estimates (\ref{expest}),
this is obviously true outside a fixed neighborhood of $\lb=0$, say
when $|\lambda| \geq 1/2$. However, since the parametrix
$N(\lb)$ is of order $1/\lb^{1/4}$ for $\lb$ close to zero, 
and the contour approaches $\lb=0$
as $\al\to 1$, we need a more detailed analysis for $|\lambda| \leq 1/2$.
In that case, we use (\ref{hout}) to write for all $\al\in[0,1)$
(in what follows the same symbols $C$ and 
$c$ stand for various positive constants independent of $\al$, $\lb$, and $n$):
\be
\left|{1\over\sqrt{\lb}}e^{-n\Re h(\lb)}\right|<
{C\over\sqrt{\mu}}\exp\left[-cn{\sqrt{\mu/ (1-\al)}\over (1+\al\mu/(1-\al))^2}\right]=
{C\over\sqrt{1-\al}}{1\over t}
\exp\left[-cn{t\over (1+\al t^2)^2}\right]\equiv f(t),
\ee
where $t=\sqrt{\mu/(1-\al)}$, $\mu=\Re\lb$. 
We need to find the maximum value of $f(t)$ in the interval
\[
t_1\equiv \sqrt{(\ep_3/\ep_2)\sin\ga_0}\le t \le {1\over\sqrt{2(1-\al)}}\equiv t_2
\]
for all $\al\in [0,1-s_0/(2n)^{2/3}]$.
For this purpose, it is convenient to consider the following two cases separately.

\noindent 1) $\al t^2 \le 1$. Then $1+\al t^2\le 2$, and we have 
\be
f(t) < {C\over t\sqrt{1-\al}}
\exp[-cnt]\equiv f_1(t).
\ee
The derivative $f'_1(t)<0$ for $t>0$, which implies
\be
\max_{t\in[t_1,1/\sqrt{\al}]}  
f(t)< f_1(t_1)< C{n^{1/3}\over s_0^{1/2}}e^{-cn}<Ce^{-cn}.\la{f1}
\ee

If $1/\sqrt{\al}>t_2$ this is all we need. Otherwise consider

\noindent 2) $\al t^2 > 1$. Then $1+\al t^2< 2\al t^2$, and we have 
\be
f(t) < {C\over t\sqrt{1-\al}}
\exp[-cn/t^3]\equiv f_2(t).
\ee
The only maximum of $f_2(t)$ is at the point $t_c=(3cn)^{1/3}$. 
Now choose sufficiently large $s_0>0$ (depending on
$\ep_2$, $\ep_3$). Then
\[
t_2={1\over\sqrt{2(1-\al)}}<cn^{1/3}/s_0^{1/2}<t_c.
\]
Therefore 
\be
\max_{t\in[1/\sqrt{\al},t_2]} f(t)< f_2(t_2)< Ce^{-cn(1-\al)^{3/2}}=
Ce^{-c\rho}.\la{f2}
\ee
Combining (\ref{f1},\ref{f2}), we finally obtain that the 
jump matrix on $(\Si_1\cup\Si_3)\setminus(\overline{U_0\cup U_1})$
is the identity up to an error of order
\be\la{outest}
\left|{1\over\sqrt{\lb}}e^{-n\Re h(\lb)}\right|<
Ce^{-c\rho}
\ee
for all $\al\in [0,1-s_0/(2n)^{2/3}]$, $2n>s_0^{3/2}$.

This estimate can be readily extended to 
complex $\al\in\cD_{\ep_0}(0)$.
The jump matrices on $\partial U_{0,1}$ admit the uniform expansions 
given by (\ref{j0},\ref{j1}).

A consequence of the above considerations is the following result: 

{\lemma \la{lemma1} Let $\rho=n|1-\al|^{3/2}$, 
$\al\in\cD_{\ep_0}(0)\cup [0,1-s_0/(2n)^{2/3}]$,
 $U=U_0\cup U_1$,
$\wt\Si=\Si^\mathrm{out}\cup \partial U$. Also let $\wt U_1$ be the circle centered at
$\lb=1$ of radius $\ep/2$. Then, for sufficiently small $\ep $,
$\ep_{j}$, $j = 0,2,3$ ($\ep_{j}$, $ j =0,2,3$ are the $\ep$-parameters
introduced above in the definition of the contour $\wt\Si$), there exists $s_{0} > 0$
such that for all $\al\in\cD_{\ep_0}(0)\cup [0,1-s_0/(2n)^{2/3}]$,
and $n > s^{3/2}_{0}/2$, 
a (unique) solution $R(\lb)$ of the $R$-RH problem exists. Moreover,
the function $R(\lb)$ admits the following asymptotic
expansion, which (and the derivative of which) is uniform 
for $\al\in\cD_{\ep_0}(0)\cup [0,1-s_0/(2n)^{2/3}]$
and all $\lb\in\wt U_1$, as $\rho\to\infty$:
\be\la{Ras}\eqalign{
R(\lb)=I+R_1(\lb)+R_2(\lb)+\cdots+R_k(\lb)+R_r^{(k+1)}(\lambda),\\
R_r^{(k+1)}(\lambda)=O(\rho^{-k-1}),\qquad 
{d\over d\lb}R_r^{(k+1)}(\lambda)=O(\rho^{-k-1}),}
\ee
$k=1,2,\dots$. The functions $R_j(\lb)=O(\rho^{-j})$
are constructed by induction as follows:
\beqa
R_1(\lb)={1\over 2\pi i}\int_{\partial U}\De_1(s){ds\over s-\lb},\qquad
R_2(\lb)={1\over 2\pi i}\int_{\partial U}(R_{1\,-}(s)\De_1(s)+\De_2(s)){ds\over s-\lb},\\
\dots,\qquad
R_k(\lb)={1\over 2\pi i}\int_{\partial U}\sum_{j=1}^k R_{k-j,-}(s)\De_j(s)
{ds\over s-\lb},\qquad R_0\equiv I.
\eeqa
}

{\bf Remark.} The uniformity means that for sufficiently small 
$\ep$, $\ep_j$, $j=0,2,3$, there exist positive constants
$s_0$, $c_1$, and $c_2$ independent of $\al$, $n$, $\lb$ such that
\be\eqalign{
|R_r^{(k+1)}|\le {c_1\over \rho^{k+1}},\qquad
\left|{d\over d\lb}R_r^{(k+1)}\right| \le {c_2\over \rho^{k+1}}
\\
\forall \lb\in\wt U_1,\quad \forall\al\in\cD_{\ep_0}(0)\cup [0,1-s_0/(2n)^{2/3}],
\quad \forall n>s_0^{3/2}/2.}
\ee
We also note that, $\rho>s_0^{3/2}/2$, $\forall\al\in\cD_{\ep_0}(0)\cup [0,1-s_0/(2n)^{2/3}]$,
$n>s_0^{3/2}/2$. 

{\bf Proof of Lemma 1.} We shall follow a similar line of argument to the one 
which was used to prove similar statement in \ci{DIKZ} (Lemma 1).
For simplicity, as in \ci{DIKZ}, we will only prove the expansion
(\ref{Ras}) in the case $k=2$, which is all that is needed for the problem at hand.
We shall also adopt the notation:
$$
R^{(3)}_{r}(\lambda) \equiv R_{r}(\lambda).
$$
Besides, as before, the symbol $c$ will stand for various positive constants independent
of $\al$, $\lb$, and $n$.

Write the jump condition for $R(\lb)$ in the form
\be\la{jlemma}
R_{0\,+}+R_{1\,+}+R_{2\,+}+R_{r\,+}=(R_{0\,-}+R_{1\,-}+R_{2\,-}+R_{r\,-})
(I+\De_1+\De_2+\De_r).
\ee
Here $\De_1$ and $\De_2$ are given by (\ref{De0},\ref{De1}) on 
$\partial U_0$, $\partial U_{1}$, respectively,
and we set $\De_1=\De_2=0$ on the rest of the contour. 
A direct analysis of the expressions (\ref{De0},\ref{De1}) shows that
$\De_k=O((n^{-k}|1-\al|^{-1/2})$ on $\partial U_0$, and 
$\De_k=O(\rho^{-k})$ on $\partial U_{1}$. 
Similarly, $\De_r=O(1/\rho^3)$ on 
$\partial U_{1}$ (this error term arises from the Bessel 
asymptotics), $\De_r=O(|1-\al|^4/\rho^3)$ on $\partial U_0$,
and, by (\ref{outest}),
$\De_r=O(e^{-c\rho})$ on  $\wt\Si\setminus \partial U$.

We now show that we can define $R_1$ and $R_2$ so that they are of order 
$1/\rho$ and $1/\rho^2$, 
respectively. We then show that the remainder $R_r$ is of order $1/\rho^3$.
Set 
\[
R_0=I.
\]
We define $R_j$ by collecting in 
(\ref{jlemma}) the terms that we want to be of the same order. First,
\be\la{jlemma1}
R_{1\,+}(\lb)=R_{1\,-}(\lb)+\De_1(\lb),\qquad \lb\in\wt\Si.
\ee
We are looking for a function $R_1(\lb)$, which is analytic outside 
$\wt\Si$, satisfying
$R_1(\lb)=O(1/\lb)$, $\lb\to\infty$, and the above jump condition.
The solution to this RH-problem is given by the Sokhotsky-Plemelj formula,
\be
R_1(\lb)=C(\De_1),
\ee
where
\[
C(f)={1\over 2\pi i}\int_{\wt\Si}f(s){ds\over s-\lb}
\]
is the Cauchy operator on $\wt\Si$.
The condition $\De_1(\lb)=O(1/\rho)$, $\lb\in\wt\Si$, 
$\rho\to\infty$ (uniform in $\al$),
implies that there exist $c,\delta,s_0>0$ such that 
\be\la{R11}
|R_1(\lb)|\le c/\rho, \quad n \geq \frac{s^{3/2}_{0}}{2}
\ee
uniformly
in $\al\in\cD_{\ep_0}(0)\cup [0,1-s_0/(2n)^{2/3}]$ and 
$\lb$ satisfying $\rm{dist}(\lb,\wt\Si)\ge \delta$.
Actually, this estimate is uniform for all $\lb\in\complex\setminus\wt\Si$ 
up to $\wt\Si$. Indeed, since
\be\la{march1}
R_{1}(\lb)={1\over 2\pi i}\int_{\partial U}\Delta_{1}(s){ds\over s-\lb},
\ee
for $\lb$ outside a fixed neighborhood of zero, this is seen by shifting the contour
to a fixed distance from the point $\lb$. Inside that neighborhood, 
the distance of the shift will depend on $\al$. Namely, the distance 
is $\ep'|1-\al|$ for a fixed (sufficiently small)  $\ep'>0$.
Then 
\be
|C(\De_1)|\le \max|\De_1| {1\over c|1-\al|}+{c\over\rho}\le
{c\over n|1-\al|^{3/2}}+{c\over\rho}
={c\over\rho}, 
\ee
on and close to $\partial U_0$. Here we used the estimate
$\Delta_{1} = O\bigl(n^{-1}\lambda^{-1/2}\bigr)$, so that
in the neighborhood of the circle $\partial U_{0}$ the inequality
$$
\max |\De_1| \leq \frac{c}{n|\sqrt{1-\al}|}
$$
holds. It should be observed
that, by the same deformation of the contour of
integration in (\ref{march1}), one obtains the analytic continuations of
both the functions $R_{1+}(\lambda)$ and $R_{1-}(\lambda)$ in the
neighborhood of the contour $\partial U$ and hence in the neighborhood of
${\wt\Si}$ (we note that on the part $\Si^\mathrm{out}$ of the
contour  ${\wt\Si}$  $R_{1}(\lambda)$ has no jump).
Moreover, the estimate
(\ref{R11}) is preserved under this analytic continuation.

Now define $R_2(\lb)$ by the jump condition
\be\la{jlemma2}
R_{2\,+}(\lb)=R_{2\,-}(\lb)+R_{1\,-}(\lb)\De_1(\lb)+\De_2(\lb),\qquad \lb\in\wt\Si,
\ee
together with the 
requirement of analyticity for $\lb\in\complex\setminus\wt\Si$, and 
the condition $R_2(\lb)=O(1/\lb)$ for $\lb\to\infty$.
The solution to this RHP is 
\be
R_2(\lb)=C(R_{1\,-}\De_1+\De_2),\qquad \lb\in\complex\setminus\wt\Si.
\ee
Using (\ref{R11}) and the estimates for $\De_2$, we obtain in the same 
way as for $R_1$, 
\be
|R_2(\lb)|\le c/\rho^2,\qquad \lb\in \complex\setminus\wt\Si,\qquad 
n \geq \frac{s^{3/2}_{0}}{2}
\ee
with the same uniformity and analyticity properties in $\al$ and $\lb$.

Now from (\ref{jlemma},\ref{jlemma1},\ref{jlemma2}) we obtain 
\be\la{jr}
R_{r\,+}(\lb)=R_{r\,-}(\lb)+M(\lb)+R_{r\,-}(\lb)\De(\lb),\qquad \lb\in\wt\Si,
\ee
where 
\[
M\equiv R_{2\,-}\De_1+(R_{1\,-}+R_{2\,-})\De_2+(I+R_{1\,-}+R_{2\,-})\De_r,\qquad
\De\equiv\De_1+\De_2+\De_r.
\]

{\bf Remark.} In the terminology of \ci{DZsobolev}, equation (\ref{jr}) is an 
inhomogeneous RH-problem of type 2. 

Since $R_r=R-I-R_1-R_2$, the matrix function $R_{r}(\lambda)$
 is analytic outside $\wt\Si$ and satisfies the condition
$R_r(\lb)=O(1/\lb)$ as $\lb\to\infty$. Therefore
\be\la{Rr}
R_r(\lb)=C(M)+C(R_{r\,-}\De),\qquad \lb\in \complex\setminus\wt\Si.
\ee
Hence
\be\la{Rr-}
R_{r\,-}(\lb)=C_-(M)+C_-(R_{r\,-}\De),\qquad \lb\in\wt\Si,
\ee
where $C_-(f)=\lim_{\lb'\to \lb}C(f)$, as $\lb'$ approaches a point $\lb\in\wt\Si$
from the $-$ side of $\wt\Si$. Now defining the operator
\[
C_{\De}(f)\equiv C_-(f\De),
\]
we represent (\ref{Rr-}) in the form
\be\la{RDe}
(I-C_{\De})(R_{r\,-})=C_-(M).
\ee
By virtue of the estimates (\ref{j0}), (\ref{j1}), 
and (\ref{outest}) we have that
\begin{equation}\label{deltal2}
||\Delta||_{L^2(\tilde{\Sigma})\cap L^{\infty}(\tilde{\Sigma})} \leq \frac{c}{\rho},
\end{equation}
for all $\al\in\cD_{\ep_0}(0)\cup [0,1-s_0/(2n)^{2/3}]$ and
$n>s_0^{3/2}/2$. 

The Cauchy operator $C_{-}$ is
bounded in the space $L^{2}(\tilde{\Sigma})$ (see, e.g., \ci{LS}), 
and by a standard scaling argument (the Cauchy operator is homogeneous of degree 0),
its norm is bounded by a constant independent of
$\alpha$. This together with the $L^{\infty}$ part
of the estimate (\ref{deltal2}) implies that the operator norm 
$||C_\De||_{L^2}=O(1/\rho)$, 
and hence $I-C_{\De}$ is invertible by a 
Neumann series for $s_0$ (and, therefore, $\rho$) sufficiently large. 
Thus (\ref{RDe}) gives
\be\la{march2}
R_{r\,-}=(I-C_{\De})^{-1}(C_-(M)),
\ee
and this proves the solvability of the $R$-RH problem for all
$\al\in\cD_{\ep_0}(0)\cup [0,1-s_0/(2n)^{2/3}]$ and
$n>s_0^{3/2}/2$. Moreover, using the $L^{2}$ part of the estimate (\ref{deltal2}),
we conclude that $\|C_{-}(M)\|_{L^2(\tilde{\Sigma})} = O(\rho^{-3})$.
Together with (\ref{march2}) this yields the uniform estimate
\be\la{march3}
\|R_{r\,-}\|_{L^2(\tilde{\Sigma})}\leq \frac{c}{\rho^3},
\ee
$$
\forall \al\in\cD_{\ep_0}(0)\cup [0,1-s_0/(2n)^{2/3}], \quad
n>s_0^{3/2}/2.
$$
The solution $R(\lambda)$ of the $R$-RH problem is given by the integral
representation
\be\la{march4}
R(\lambda) = I + R_{1}(\lambda) + R_{2}(\lambda) + C(M) 
+ C(R_{r\,-}\Delta)(\lambda),
\ee
$$
\lambda \in \mathbb C \setminus \tilde{\Sigma}.
$$

{\bf Remark.} Let $\Omega_{k}$, $ k = 1,2,3,4$ denote
the connected components of the set ${\mathbb C} \setminus \tilde{\Sigma}$.
Then, using again the possibility of
the contour deformation when solving the integral equation (\ref{RDe}),
and taking into account the triviality of the jump matrix monodromy at each
node point of the contour $\tilde{\Sigma}$,  
we conclude that the restriction $R|_{\Omega_{k}}(\lambda)$ is
continuous in $\overline{\Omega_{k}}$ for each $k$ (see e.g. \cite{BDT}). 
This means that equation (\ref{march4}) defines the solution of
the $R$-RH problem in the classical, point-wise continuous, sense.

Combining the inequality (\ref{march3}) with equation (\ref{march4}),
we can complete the proof of the lemma. Indeed, assuming
that $\lambda \in \wt U_1$, we immediately obtain the estimate
\be\la{march5}
|C(M)(\lambda)| \leq \frac{c}{\rho^3}, \quad  n > s_0^{3/2}/2,
\ee
for the fourth term in the r.h.s. of (\ref{march4}), and the estimate
\be\la{march6}
|C(R_{r\,-}\Delta)(\lambda)| \leq c
||R_{r\,-}||_{L^2(\tilde{\Sigma})}||\Delta||_{L^2(\tilde{\Sigma})}\leq
\frac{c}{\rho^3},
\ee
$$
\quad n > s_0^{3/2}/2,
$$
for the fifth term. Both the estimates are uniform in
$\al\in\cD_{\ep_0}(0)\cup [0,1-s_0/(2n)^{2/3}]$. Together they
yield the estimate
\be\la{march66}
|R_{r}(\lambda)| \leq \frac{c}{\rho^3}, \quad n > s_0^{3/2}/2,
\ee
uniformly in $\al\in\cD_{\ep_0}(0)\cup [0,1-s_0/(2n)^{2/3}]$ and 
$\lb$ lying in $\wt U_1$. This establishes part of the estimate (\ref{Ras})
for the error term. The estimate for the derivative follows immediately
from (\ref{Rr}).
This completes the proof of the lemma (in the case $k=2$). 
$\Box$
\vskip .2in

{\bf Remark.} (Cf. Remark 2 in \ci{DIKZ}.)  
Part of the assertion of Lemma 1 is that the solution of the $R$-RH problem,
and hence of the original $T$-RH problem,
exists and is unique for all $\al\in\cD_{\ep_0}(0)\cup [0,1-s_0/(2n)^{2/3}]$ and 
$n > s_0^{3/2}/2$ with $s_{0}$ sufficiently large.
This is all we need in the analysis that follows;
however, the solution of the $R$-RH problem actually exists and is unique for
all $\al\in\cD_{\ep_0}(0)\cup [0,1)$ and {\it all}  $n > 0$
for some (possibly smaller) $\ep_0 > 0$. Indeed, by the discussion following
(\ref{TT}), the $T$-RH problem, and hence the $R$-RH problem, is solvable
for all $\alpha \in [0, 1)$, $n > 0$. 
Since, by the previous remark, the solution of the $R$-RH problem 
is continuous up to the contour, the problem is easily seen to be solvable for  
$\al\in\cD_{\ep'_0}(0)$, $0 < n \le s_0^{3/2}/2$ for some $\ep'_0 > 0$
by continuity of the jump matrix at $\alpha = 0$. By Lemma 1,
the $R$-RH problem is solvable for all $\al\in\cD_{\ep_0}(0)$,
$n > s_0^{3/2}/2$. Thus the  $R$-RH problem,
and hence the $T$-RH problem, is solvable for all $ n > 0$ on
$\cD_{\ep''_0}(0)\cup [0,1)$, where $\ep'' = \min\{\ep_0,\ep'_{0}\}$.


\section{Evaluation of the differential identity}
\subsection{Exact transformations}
We start with the differential identity (\ref{diffV}).
Note that since $V(z)$ is related to $U(\lb)$ by the expression (\ref{U},\ref{zlb})
\[
U(\lb)=V(z(\lb)),\qquad z={\al\lb\over 1-\al+\al\lb},
\]
we have
\be
\left.{d\lb\over dz}\right|_{z=\al}={1\over\al(1-\al)},
\ee  
and (\ref{diffV}) can be rewritten in terms of $U(\lb)$ as follows
\be\la{diffU}
{d\over d\al}\ln D_n(\al)={e^{-4n\al}\over 2\pi i\al(1-\al)}
(U_{11}(1)U'_{21}(1)-U'_{11}(1)U_{21}(1)).
\ee
Note that the derivatives in (\ref{diffU}) are taken w.r.t. $\lb$.

By (\ref{T},\ref{TT}), 
the matrix elements of $U(\lb)$ can be expressed in terms of $T(\lb)$ as follows:
\be\eqalign{
U_{11}(\lb)=\left[T^{-1}\left(-{1-\al\over \al}\right)T(\lb)\right]_{11}e^{n\hat g(\lb)},
\\
U_{21}(\lb)=\left[T^{-1}\left(-{1-\al\over \al}\right)T(\lb)\right]_{21}e^{-nl}
e^{n\hat g(\lb)}}.
\ee
Furthermore, for $\lb$ outside the lens in $U_1$
\be
T(\lb)=S(\lb),\qquad S(\lb)=R(\lb)P_1(\lb). 
\ee
Note also that by (\ref{P1p})
\[
S_{j1}=(R(\lb)E_n(\lb)\Psi(\ze))_{j1}e^{n\phi(\lb)},\qquad j=1,2,
\]
and, as follows from the definitions of the functions $\phi$, $h$, and the properties of $g(z)$,
\[
\phi(\lb)+\hat g(\lb)=\mp {1\over 2}h+\hat g=
\mp{\hat g_+ -\hat g_-\over 2}+ \hat g_{\pm}=
 {\hat g_+ +\hat g_-\over 2}=2z(\lb)+l/2,
\]
where $\hat g_{\pm}(\lb)$ stand for the analytic continuation of these functions.
Here the upper sign corresponds to $\Im\lb>0$, and the lower, to $\Im\lb<0$.

Hence, (\ref{diffU}) finally gives
\be\la{diffR}
{d\over d\al}\ln D_n(\al)={1\over 2\pi i\al(1-\al)}
((RE\Psi)_{11}(1)(RE\Psi)'_{21}(1)-(RE\Psi)'_{11}(1)(RE\Psi)_{21}(1)),
\ee
where we used the fact that $\det T^{-1}(-(1-\al)/\al)=1$. 
In (\ref{diffR}), the derivative at $\lb=1$ is taken along a path in $U_1$ outside the lens.
In the next subsection we use the solution of the Riemann-Hilbert problem for $R(\lb)$ 
(found in Section 3) to construct the asymptotics of the r.h.s. of (\ref{diffR}).

\subsection{Asymptotics}
Consecutive asymptotic terms in the expansion of the logarithmic derivative (\ref{diffR})
are generated by consecutive terms in (\ref{Ras}):
\[
R(\lb)=I+R_1(\lb)+R_2(\lb)+\cdots.
\]
Thus, setting $R=I$ in (\ref{diffR}) gives the main asymptotic term of ${d\over d\al}\ln D_n(\al)$:
\be\la{1termp}
{1\over 2\pi i\al(1-\al)}
((E\Psi)_{11}(1)(E\Psi)'_{21}(1)-(E\Psi)'_{11}(1)(E\Psi)_{21}(1)),
\ee
Using (\ref{E1}) and (\ref{Q3}), we obtain
\be
(E\Psi)_{11}(\ze)=\mu_+(\lb),\qquad (E\Psi)_{21}(\ze)=-i\mu_-(\lb),
\ee
where
\be
\mu_{\pm}(\lb)=
\sqrt{\pi\over 2}\ze^{1/4}(m^{-1}(\lb)I_0(\sqrt{\ze})\pm m(\lb)I'_0(\sqrt{\ze})).
\ee
Using the expansion of Bessel functions as $\ze\to 0$ (i.e. $\lb\to 1$), we obtain
\be\eqalign{
\mu_{\pm}(1)\equiv M=\sqrt{\pi n}(1-\al)^{3/4},\qquad
\mu'_\pm (1)\equiv a\pm b,\\ 
a=M \left[n^2(1-\al)^3-{\al\over 2}+{1\over 6}\right],\qquad
b=M n(1-\al)^{3/2}.}\la{mu}
\ee
Substituting these values into (\ref{1termp}), we find the main asymptotic term
\be
{d\over d\al}\ln D_n(\al)\sim {n^2\over\al}(1-\al)^2.
\ee
To obtain the next term, we need to compute first
\be\la{R1l}
R_1(1)={1\over 2\pi i}\int_{\partial U}{\De_1(\lb)\over \lb-1}d\lb,\qquad
R'_1(1)={1\over 2\pi i}\int_{\partial U}{\De_1(\lb)\over (\lb-1)^2}d\lb.
\ee
We now examine $\De_{1,2}$ in the neighborhoods
of the points $\lb=0,1$.
Using (\ref{ze0}) and expanding the matrix elements of $N(\lb)$, 
we obtain from (\ref{De0}):
\be\eqalign{
\De_1(\lb)={C_1\over\lb}+{\sqrt{1-\al}\over 32n (1+\al)}\pmatrix{
F_0(\al)-5/2 & i(F_0(\al)+7/2)\cr i(F_0(\al)+7/2)& -F_0(\al)+5/2}+O(\lb),\\
F_0(\al)={1-6\al-3\al^2\over 6(1-\al^2)},\qquad
C_1={\sqrt{1-\al}\over 32n (1+\al)}\pmatrix{-1& -i\cr -i & 1},\qquad \lb\in U_0.}\la{C1}
\ee

For $\De_2(\lb)$, we obtain similarly:
\be
\De_2(\lb)={3(1-\al)(1+O(\lb))\over 2^9n^2(1+\al)^2\lb}
\pmatrix{1 & 4i\cr -4i & 1},\qquad
\lb\in U_0.
\ee

In $U_1$, a similar calculation based on (\ref{De1}) and (\ref{ze1}) gives
($\lb=1+u$)
\be\eqalign{
\De_1(\lb)={A_1\over u}+{1\over 32n(1-\al)^{3/2}}
\pmatrix{-5/2+\al+1/6 & -i(7/2+\al+1/6)\cr
-i(7/2+\al+1/6) & -(-5/2+\al+1/6)}+\\
{u\over 32n(1-\al)^{3/2}}\pmatrix{
3/2-(5/2)(\al+1/6)+F_1(\al) & -i(-3/2+(7/2)(\al+1/6)+F_1(\al))\cr
-i(-3/2+(7/2)(\al+1/6)+F_1(\al)) & -(3/2-(5/2)(\al+1/6)+F_1(\al))}+\\
O(u^2),\qquad 1+u=\lb,\qquad \lb\in U_1,\qquad
A_1={1\over 32n(1-\al)^{3/2}}\pmatrix{1 & -i \cr -i & -1},}\la{A1}
\ee
\be\la{De2as1}
F_1(\al)={\al\over 6}-{31\over 4\cdot 45},\qquad
\De_2(\lb)={3(1+(2\al+1/3)u+O(u^2))\over 2^9 n^2(1-\al)^3 u} 
\pmatrix{-1 & 4i\cr -4i & -1},\qquad
\lb\in U_1.
\ee

Now the expressions for $R(1)$ and $R'(1)$ are obtained from the above results 
and (\ref{R1l}) by a straightforward residue calculation:
\be\eqalign{
R_1(1)=\pmatrix{\de & \eta \cr \eta & -\de},\qquad
\de={1\over 32n(1-\al)^{3/2}}\left[5/2-\al-1/6-{(1-\al)^2\over 1+\al}\right],\\
\eta={i\over 32n(1-\al)^{3/2}}\left[7/2+\al+1/6-{(1-\al)^2\over 1+\al}\right],\\
R'_1(1)=\pmatrix{\si & \tau \cr \tau & -\si},\qquad
\si={1\over 32n(1-\al)^{3/2}}\left[-3/2+(5/2)(\al+1/6)-F_1(\al)+
{(1-\al)^2\over 1+\al}\right],\\
\tau={i\over 32n(1-\al)^{3/2}}\left[-3/2+(7/2)(\al+1/6)+F_1(\al)+
{(1-\al)^2\over 1+\al}
\right].}
\ee
Note that the contours $\partial U_{0,1}$ are traversed in the negative
direction.

We shall be using the following notation for the expansion terms of 
the logarithmic derivative (\ref{diffR}).
We denote $R_k\cdot R_m$ ($R_0\equiv I$) the term given by
\be\eqalign{
{1\over 2\pi i\al(1-\al)(1+\de_{k,m})}
((R_k E\Psi)_{11}(1)(R_m E\Psi)'_{21}(1)+
(R_m E\Psi)_{11}(1)(R_k E\Psi)'_{21}(1)\\
-(R_k E\Psi)'_{11}(1)(R_m E\Psi)_{21}(1)-
(R_m E\Psi)'_{11}(1)(R_k E\Psi)_{21}(1)).}
\ee
For example, the main term (\ref{1termp}) is $I\cdot I$.
We can now evaluate the next ($R_1\cdot I$) term in the expansion. 
It is written as follows:
\be\eqalign{
{d\over d\al}\ln D_n(\al)-{n^2\over\al}(1-\al)^2\sim
{1\over 2\pi i\al(1-\al)}\left(\left\{R_1(1)\pmatrix{1\cr -i}\right\}_1 M (-i\mu'_-(1))+
\right.\\ 
M \left\{R_1(\lb)\pmatrix{\mu_+(\lb)\cr -i\mu_-(\lb)}\right\}'_2(1)-
\left\{R_1(\lb)\pmatrix{\mu_+(\lb)\cr -i\mu_-(\lb)}\right\}'_1(1)(-iM)-\\
\left.\mu'_+(1)\left\{R_1(1)\pmatrix{1\cr -i}\right\}_2 M\right)=
{M^2 (\tau+i\si)\over \pi i\al(1-\al)}={\al\over 4(1-\al^2)},}\la{2term}
\ee
where we first simplified the expression substituting the above symbolic representation 
of $R_1$ in terms of $\de$, $\eta$, $\si$, $\tau$, and used their numerical values only 
at the last step.

It turns out that the two terms in the asymptotics just obtained is all we need 
(up to the error term).
The following lemma is the main result of this section:

{\lemma\la{lemma2}
There exists $s_0>0$ such that the expansion
\begin{eqnarray}\la{Lemma}
{d\over d\al}\ln D_n(\al)={n^2\over\al}(1-\al)^2+{\al\over 4(1-\al^2)}+
r(n,\al),\la{di1}\\
r(n,\al)={1\over 1-\al}O\left({1\over\rho}\right),\qquad
\rho=n|1-\al|^{3/2},\la{di2}
\end{eqnarray}
holds uniformly in $\al\in(0,1-s_0/(2n)^{2/3}]$ for all $n>s_0^{3/2}/2$. 
}

{\bf Proof.} It only remains to prove the expression for the error term.
We consider the expansion of $R(\lb)$ up to the third term:
$R=I+R_1+R_2+R_r$.
Since $R_k=O(\rho^{-k})$ and, according to (\ref{mu}),
$\mu'_\pm(1)\mu_\pm(1)=O(\rho^3)$, it is not difficult to deduce from
(\ref{diffR}) (cf. (\ref{2term})) that the contribution of the terms 
$R_r\cdot R_1$, $R_2\cdot R_2$ and higher are of order
$(\al(1-\al))^{-1}O(\rho^{-1})$. 
Thus we shall need to consider in detail only the following 4 terms:
$R_1\cdot R_1$, $R_2\cdot I$, $R_2\cdot R_1$, $R_r\cdot I$.

For the $R_1\cdot R_1$ term, which we denote $L_{11}$, we obtain after 
a calculation similar to 
(\ref{2term}):
\be\la{L11}
L_{11}=-{n^2\over\al}(1-\al)^2(\de^2+\eta^2)=
{1\over 2^8\al(1-\al^2)}(\al+2/3)(2+5\al-\al^2).
\ee

For further analysis, we need to calculate $R_2(1)$. It is given by the formula:
\be\la{R2(1)}
R_2(1)={1\over 2\pi i}\int_{\partial U} {R_{1-}(\lb)\De_1(\lb)+\De_2(\lb)\over \lb-1}d\lb.
\ee
The solution of the Riemann-Hilbert problem for $R_1$ inside $U_{1,0}$ is given by the 
expression (which we write on the boundary)
\be\la{R1-}
R_{1-}(\lb)={A_1\over \lb-1}+{C_1\over \lb}-\De_1 (\lb),
\qquad \lb\in\partial U,
\ee
where $A_1$, $C_1$ are defined in (\ref{A1},\ref{C1}).
Note that outside $U_{1,0}$ the solution is
\[
R_{1}(\lb)={A_1\over \lb-1}+{C_1\over \lb}.
\]
It is easily seen that the jump, analyticity conditions, and the condition at infinity
of the Riemann-Hilbert problem for $R_1(\lb)$ are satisfied, and therefore, 
by uniqueness, this is the solution. 

The expansions for $\De_{1,2}$ obtained above and the formulas (\ref{R1-},\ref{R2(1)})
give, by a residue calculation, the final expression for $R_2(1)$:
\be\eqalign{
R_2(1)=\pmatrix{\ga & -\beta \cr \beta & \ga},\qquad
\ga={-1\over 2^9n^2(1-\al)^3}\left[(3\al-1)\left(1-{(1-\al)^2\over 3(1+\al)}\right)
+3-{(1-\al)^2\over 1+\al}\right],}
\ee
where the expression for $\beta$ is omitted as it is not needed below.

To compute the ``$R_2\cdot I$'' term (which we denote $L_{20}$) note first 
that the contribution of 
the terms in that expression involving $R'_2(1)$ is of order 
$(\al(1-\al))^{-1}O(\rho^{-1})$ and we need not calculate them. The remainder gives
a nontrivial contribution, and we obtain:
\be
L_{20}={2bM\ga\over \pi\al(1-\al)}+
{1\over \al(1-\al)}O\left({1\over n(1-\al)^{3/2}}\right).
\ee
The expression for $\ga$ tells us that this is equal to $-L_{11}$ (\ref{L11}) up to 
the error term. Thus, we conclude that the contributions of $R_2\cdot I$ and 
$R_1\cdot R_1$ terms cancel each other.

The analysis of the $R_2\cdot R_1$ term is now easy to carry out, and we find that this 
term is of order $(\al(1-\al))^{-1}O(\rho^{-2})$. 

For {\it any} matrix elements of $R_r(1)$ (we only know they are of order  
$O(\rho^{-3})$), we obtain that the $R_r\cdot I$ term is of order 
$(\al(1-\al))^{-1}O(\rho^{-1})$.

Thus, in view of uniformity of the error term in the expansion of $R(\lb)$, 
the lemma is proven but with the remainder
\be\la{err}
r(n,\al)={1\over \al(1-\al)}O\left({1\over\rho}\right).
\ee
We now show that $\al$ in the denominator here can be omitted. 
First, we notice that $r(n,\al)=O_n(1)$ as $\al\to 0$ and $n$ is fixed:
this follows immediately after substitution of the expansion (\ref{Dal0})
into the l.h.s. of (\ref{Lemma}). However, we need an estimate which is uniform in $n$.
To obtain such an estimate, we use the extensions of our expressions for complex $\al$
discussed above. As follows from (\ref{Lemma},\ref{Dal0}), $r(n,\al)$ is an analytic function of
$\al$ in $\cD_{\ep_0}(0)$. Thus
\be
r(n,\al)={1\over 2\pi i}\int_{\partial \cD_{\ep_0 /2}(0)}{r(n,\wt\al)\over \wt\al-\al}d\wt\al,
\qquad |\al|<\ep_0 /4.
\ee
Since by (\ref{err}),
$r(n,\wt\al)$ is uniformly bounded on $\partial\cD_{\ep_0/2}(0)$, it follows
that $r(n,\al)$ is uniformly bounded by $O(1/\rho)$ 
for all $\al\in\cD_{\ep_0 /4}(0)$, and all $n>s_0^{3/2}/2$.
Lemma 2 is proven. $\Box$

\section{Proof of Theorem 1}
Integrating the differential identity (\ref{Lemma}) from $\al_0$ (close to zero from above)
to any $\al_0<\al\le 1-s_0/(2n)^{2/3}$, we obtain:
\be\la{intD}
\ln D_n(\al)- \ln D_n(\al_0)=
n^2\left(\ln{\al\over\al_0}-2(\al-\al_0)+{\al^2-\al_0^2\over 2}\right)
-{1\over 8}\ln{1-\al^2\over 1-\al_0^2}+
O\left(1\over n(1-\al)^{3/2}\right)
\ee
for all $n>s_0^{3/2}/2$. Note from (\ref{di2}) that the term $O(1/n(1-\al)^{3/2})$ does not depend on $\al_0$. 
Substituting for $\ln D_n(\al_0)$ the expansion (\ref{Dal0}) and taking the limit 
$\al_0\to 0$, we obtain for any $0<\al\le 1-s_0/(2n)^{2/3}$,
\be\eqalign{
\ln D_n(\al)=
n^2\left({3\over 2}+\ln\al-2\al+{\al^2\over 2}\right)-{1\over 12}\ln n -
{1\over 8}\ln (1-\al^2)+\\
{1\over 12}\ln2 +\ze'(-1)+
O\left(1\over n(1-\al)^{3/2}\right)+\de_n.}\la{intD2}
\ee

Fix any $s>s_0$ and, for $n$ sufficiently large, set
$\al=1-s/(2n)^{2/3}$. Now take the limit $n\to\infty$. 
As $n\to\infty$, the r.h.s. of (\ref{intD2}) becomes
\be
-{s^3\over 12}-{1\over 8}\ln s + {1\over 24}\ln 2 +\ze'(-1)+O(s^{-3/2}).
\ee
On the other hand, as $s$ is any fixed number $s>s_0$, the l.h.s. of (\ref{intD2})
converges to $\ln\det(I-K_s)$ by (\ref{detA}).
$\Box$

\bigskip
\bigskip
 
\noindent{\Large\bf Acknowledgements}
\bigskip

\noindent
Percy Deift was supported
in part by NSF grant \# DMS 0500923 and also by a Friends of
the Institute Visiting Membership at the Institute for Advanced Study in Princeton,
Spring 2006. Alexander Its was supported
in part by NSF grant \# DMS-0401009. 
The authors thank V. Tarasov for a useful discussion which
took place after one of us (A.I.) gave a talk on our preceding work
\ci{DIKZ}. In fact, it is during this discussion that the idea
to use  the Airy-limit of the classical polynomials in order to prove 
the Tracy-Widom conjecture in the framework
of the approach of \ci{DIKZ} was born.


\section{Appendix}
Here we present an alternative derivation of the identity (\ref{diffV}).
Let 
\be
\phi(x)={1\over 2}\om_{n-1}(x),\qquad \psi(x)={1\over 2}\om_n (x).
\ee
The determinant (\ref{det}) is written then as follows:
\be
D_n(\al)=\det\left(I-{\phi(x)\psi(y)-\phi(y)\psi(x)\over x-y}
\chi_{(\al,\infty)}\right).\la{det2}
\ee

The operator $K(x,y)=(\phi(x)\psi(y)-\phi(y)\psi(x))/(x-y)$ is of integrable type, 
and hence (see, e.g., \ci{IIKS,Dintop,DIZ}) $D_n(\al)$ 
is related to the following Riemann-Hilbert problem
for a $2\times 2$ matrix-valued function $Y(z)$ (Figure \ref{1}):

\begin{figure}
\centerline{\psfig{file=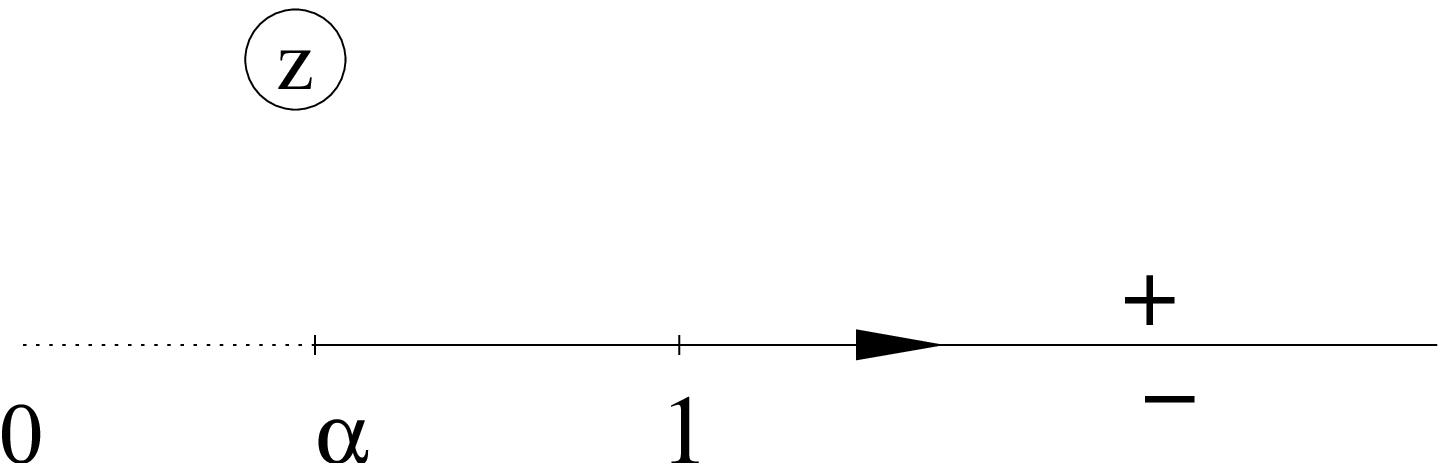,width=3.0in,angle=0}}
\vspace{0cm}
\caption{
Contour for the $Y$-Riemann-Hilbert problem.}
\label{1}
\end{figure}

\begin{enumerate}
    \item[(a)]
        $Y(z)$ is  analytic for $z\in\bbc \setminus [\al,\infty)$.
    \item[(b)] 
Let $x\in(\al,\infty)$.
$Y(z)$ has $L_2$ boundary values
$Y_{+}(x)$ as $z$ approaches $x$ from
above, and $Y_{-}(x)$, from below. 
They are related by the jump condition
\begin{equation}\label{RHPYb}
            Y_+(x) = Y_-(x) v_Y(x),\qquad
            v_Y(x)=\pmatrix{
                1+2\pi i\phi(x)\psi(x) & -2\pi i \psi(x)^2 \cr
                2\pi i \phi(x)^2 & 1-2\pi i\phi(x)\psi(x)},
            \qquad\mbox{$x\in(\al,\infty)$.}
        \end{equation}
    \item[(c)]
        $Y(z)$ has the following asymptotic behavior at infinity:
        \begin{equation} \label{RHPYc}
            Y(z) = I+ O \left( \frac{1}{z} \right),
            \qquad \mbox{as $z\to\infty$.}
        \end{equation}
\end{enumerate}

As in \ci{BD, Alexander}, it is possible to reduce the RHP for Y to 
an equivalent RHP with an ``elementary'',
in fact constant, jump matrix (see (\ref{Xjump}) below).
Note first that for any functions $\ti\psi(x)$, $\ti\phi(x)$ such that
$\psi(x)\ti\phi(x)-\phi(x)\ti\psi(x)=1$, we have
\be
v_Y(x)=A(x)\pmatrix{1 & -2\pi i\cr 0 & 1}A^{-1}(x),\qquad 
A(x)=\pmatrix{\psi(x) & \ti\psi(x) \cr \phi(x) & \ti\phi(x)},\la{j}
\ee
Note that the condition on $\ti\psi(x)$, $\ti\phi(x)$ is equivalent to the following one:
\[
\det A(x)=1.
\]

Let 
\be
\Phi(z)=\pmatrix{\psi(z) & e^{2nz}\int_0^\infty {\psi(\xi)\over \xi-z} e^{-2n\xi}d\xi\cr
 \phi(z) & e^{2nz}\int_0^\infty {\phi(\xi)\over \xi-z} e^{-2n\xi}d\xi}.
\ee
The function $\Phi(z)$ is analytic in $\bbc\setminus\bbr_+$.
Using the orthogonality property of the polynomials $p_n(x)$, $p_{n-1}(x)$ with 
respect to the weight 
$e^{-4nx}$, we see that $\Phi(z)$ solves the following RHP on $\bbr_+$:
\begin{enumerate}
    \item[(a)]
        $\Phi(z)$ is  analytic for $z\in\bbc \setminus [0,\infty)$.
    \item[(b)] 
For $x\in (0,\infty)$ the $L_2$ boundary values
$\Phi_{+}(x)$ and $\Phi_{-}(x)$ are related by the jump condition
\begin{equation}
            \Phi_+(x) = \Phi_-(x)
            \pmatrix{
                1 & 2\pi i \cr
                0 & 1},
            \qquad\mbox{$x\in (0,\infty)$.}
        \end{equation}
    \item[(c)]
        $\Phi(z)$ has the following asymptotic behavior as $z\to\infty$:
        \begin{equation}
            \Phi(z) = \left(I+ O \left( \frac{1}{z} \right)\right)
            \left( {\ka_n\over 2}e^{-2nz}z^n\right)^{\si_3}.
    \end{equation}
\end{enumerate}

By standard arguments, see \ci{Dbook}, $\det\Phi(z)=1$. Hence, we see that for $x>0$,
we can take in (\ref{j})
\be
A(x)=\Phi_+(x).
\ee

The decomposition (\ref{j}) suggests the following transformation of the 
Riemann-Hilbert problem.
Let
\be\la{X}
X(z)=Y(z)\Phi(z),
\ee

It is easy to verify that $X(z)$ satisfies the following problem:
\begin{enumerate}
    \item[(a)]
        $X(z)$ is  analytic for $z\in\bbc \setminus [0,\al]$.
    \item[(b)] 
For $x\in (0,\al)$ the $L_2$ boundary values
$X_{+}(x)$ and $X_{-}(x)$ are related by the jump condition
\begin{equation}\la{Xjump}
            X_+(x) = X_-(x)
            \pmatrix{
                1 & 2\pi i \cr
                0 & 1},
            \qquad\mbox{$x\in (0,\al)$.}
        \end{equation}
    \item[(c)]
        $X(z)$ has the following asymptotic behavior as $z\to\infty$:
        \begin{equation} 
            X(z) = \left(I+ O \left( \frac{1}{z} \right)\right)
            \left( {\ka_n\over 2}e^{-2nz}z^n\right)^{\si_3}, 
    \end{equation}
\end{enumerate}
Thus $X(z)$ satisfies the same RHP as $\Phi(z)$, but now on the interval $(0,\al)$.

The transformation 
\be\la{V}
V(z)=\left(\sqrt{2\pi i}{\ka_n\over 2}\right)^{-\si_3} X(z) e^{2nz\si_3}
\left(2\pi i\right)^{\si_3/2}
\ee
converts the RHP to the RHP for $V(z)$ of Section 3.

We now turn to the derivation of the identity for $D_n(\al)$.
Write the determinant (\ref{det2}) in the form
\[
D_n(\al)=\det(I-K),
\]
where $K$ is an integral operator acting on functions $f(x)$ from $L^2(\al,\infty)$
as follows:
\[
(Kf)(x)=\int_\al^{\infty}K(x,y)f(y)dy,\qquad K(x,y)={\phi(x)\psi(y)-\phi(y)\psi(x)
\over x-y}.
\]
The logarithmic derivative of $D_n(\al)$ w.r.t. $\al$ has the form
\be\eqalign{
{d\over d\al}\ln D_n(\al)=-\mathrm{tr}\, \left((I-K)^{-1}{dK\over d\al}\right)=
((I-K)^{-1}K)(\al,\al)=\\ 
((I-K)^{-1}(K-I+I))(\al,\al)=R(\al,\al),}\la{diff1}
\ee
where $R(x,y)$ is the kernel of the operator $(I-K)^{-1}-I$.
As noted above, the kernel $K(x,y)$ has the structure of an ``integrable'' kernel. 
A consequence of this fact is the identity
\be
R(x,y)={-F_1(x)F_2(y)+F_2(x)F_1(y)\over x-y},
\ee
where the $F_j(z)$ are expressed in terms of the solution of the Riemann-Hilbert problem
for $Y(z)$ as follows:
\be
F_j(z)=Y_{+,j1}\psi+Y_{+,j2}\phi,\qquad j=1,2.
\ee
Comparing this with the definition (\ref{X}) of $X(z)$ we see that
\be
F_j(z)=X_{j1}(z),\qquad j=1,2.
\ee
Substituting then $R(\al,\al)=\lim_{x\to\al}R(x,\al)$ into (\ref{diff1}),
we obtain:
\be\la{diffX}
{d\over d\al}\ln D_n(\al)=X_{11}(\al)X'_{21}(\al)-X'_{11}(\al)X_{21}(\al),
\ee
which expresses the logarithmic derivative of $D_n(\al)$ 
in terms of the solution of the Riemann-Hilbert problem for $X(z)$. 
Now the function $X(z)$ is related to $V(z)$ by the expression (\ref{V}). In particular,
\[
X_{11}(z)={\ka_n\over 2}e^{-2nz}V_{11}(z),\qquad
X_{21}(z)={1\over \pi i\ka_n}e^{-2nz}V_{21}(z).
\]
Calculating the derivatives of these quantities at $z=\al$ and substituting into 
(\ref{diffX}), we finally obtain (\ref{diffV}).

\end{document}